\documentclass[12pt]{article}
\usepackage{amssymb}
\usepackage{amsmath}
\usepackage{amsopn}
\numberwithin{equation}{section}

\newtheorem{theorem}{Theorem}
\newtheorem{proposition}{Proposition}[section] 
\newtheorem{corollary}[proposition]{Corollary}
\newtheorem{lemma}[proposition]{Lemma}

\let\ds=\displaystyle
\def\coun{\varepsilon}
\def\bn{\begin{equation}}
\def\ed{\end{equation}}

\def\<{\langle}
\def\>{\rangle}
\newcommand{\CC}{{\mathbb C}}

\newcommand{\ZZ}{{\mathbb Z}}

\def\r#1{(\ref{#1})}
\let\rf=\r
\def\ot{\otimes}
 \def\d{\delta}

\def\sk#1{\left(#1\right)}
\def\U{\overline{U}}

\let\nn=\nonumber

\def\Pfp{{P}^+}
\def\Pfm{{P}^-}
\def\hPfp{{\hat P}^+}
\def\hPfm{{\hat P}^-}

\def\Uqgln{U_q(\widehat{\mathfrak{gl}}_N)}
\def\Uqgl#1{U_q(\widehat{\mathfrak{gl}}_{#1})}
\def\Uqbpp#1{U_q(\widehat{\mathfrak{b}}^+_{#1})}

\def\Uqsl2{U_q(\widehat{\mathfrak{sl}}_2)}

\def\ee{{E}}\def\ff{{F}}

\def\EE{{\rm E}}\def\FF{{\rm F}}

\def\prodl{\mathop{\overleftarrow\prod}\limits}

\def\ct{\mathbb{T}}
\def\End{\textrm{End}}

\def\bbb{\mathbb{B}}

\def\ccR{\mathbb{R}}
\def\si{\sigma}
\def\RR{{\rm R}}
\def\LL{{\rm L}}
\def\ord{\prec}
\def\E{{\rm E}}
\def\F{{\mathcal{F}}}
\def\SF{{\mathcal{S}}}

\def\qsym{\varpi}
\def\tSym{\overline{\rm Sym}}
\def\Sym{{\rm Sym}}
\def\ele{{\sf E}}

\def\Uqqq{U_q(\mathfrak{gl}_N)}
\def\weig{\Lambda}

\def\lll{l}
\def\rr{r}
\def\ss{s}

\def\admis#1{[[\bar #1]]}
\def\seg#1#2{[#2,#1]}
\def\segg#1{[#1]}
\def\meg#1#2{[#1,#2]}
\def\megg#1{[#1]}
\def\bms#1{\bar\ms_{#1}}

\def\ppp{\tilde{\mathbf{\ss}}}
\def\Ff{f}
\def\acc{G}
\let\qsym=\beta

\def\qSym{\overline{\rm Sym}^{\ (q)}}

\def\ms{{m}}

\def\tFF{\hat{\rm F}}
\def\tk{\hat{k}}
\def\tEE{\hat{\rm E}}
\def\tiF{\hat{F}}
\def\tiE{\hat{E}}
\def\ticalW{\hat{\mathcal{W}}}
\def\calW{{\mathcal{W}}}

\def\ticSF{\hat{\mathcal{S}}}

\def\hU{\hat{U}}
\def\bm{\bar{\mathbf{m}}}
\def\bs{\bar{\mathbf{s}}}
\def\bfm{{\mathbf{m}}}
\def\bfs{{\mathbf{s}}}
\def\hF{\hat{\mathcal{F}}}
\def\hZ{\hat{Z}}
\def\hS{\hat{S}}
\def\hI{\hat{I}}
\def\qed{\hfill$\square$\medskip}

\begin{document}

\begin{center}

\hfill ITEP-TH-55/07\\
\bigskip
{\Large\bf On the universal weight function for the quantum\\ affine
algebra $\Uqgln$}
\par\bigskip\medskip
{\bf
Andrey Oskin$^{\bullet}$\footnote{E-mail: aoskin@theor.jinr.ru},
Stanislav Pakuliak$^{\star\bullet}$\footnote{E-mail: pakuliak@theor.jinr.ru}
and Alexey Silantyev$^{\bullet\circ}$\footnote{E-mail: silant@tonton.univ-angers.fr}
}\ \
\par\bigskip\medskip
$^\star${\it Institute of Theoretical \& Experimental Physics, 117259
Moscow, Russia}
\par\smallskip
$^\bullet${\it Laboratory of Theoretical Physics, JINR,
141980 Dubna, Moscow reg., Russia}\\
$^\circ${\it D\'epartment de Math\'ematiques,
Universit\'e d'Angers,\\
2 Bd. Lavoisier, 49045 Angers, France}\\
\bigskip

\end{center}

\thispagestyle{empty}

\begin{abstract}
We continue investigation of the universal weight function for the
quantum
affine algebra  $U_q(\widehat{\mathfrak{gl}}_N)$ started in \cite{KPT} and \cite{KP-GLN}.
We obtain two recurrence relations for the universal weight function
applying the method of projections developed in \cite{EKhP}.
On the level of the evaluation representation of $U_q(\widehat{\mathfrak{gl}}_N)$
we reproduce both recurrence relations for the off-shell Bethe vectors calculated
in \cite{TV3} using combinatorial methods.
\end{abstract}


\setcounter{footnote}{0}

\section{Introduction}

Hierarchical (nested) Bethe ansatz was designed in \cite{KR83} to
construct the eigenvectors of the commuting integrals for quantum
integrable models associated with the Lie algebra
$\mathfrak{gl}_N$. It is based on the inductive procedure which
relates $\mathfrak{gl}_N$ and $\mathfrak{gl}_{N-1}$ Bethe vectors.
If the parameters of these, so called off-shell Bethe vectors,
satisfy Bethe equations, then the corresponding vectors are
eigenvectors of the commuting set of operators in some quantum
integrable model.

Further development of the  off-shell Bethe vectors theory was
achieved in \cite{VT}, where they were presented as particular
matrix elements of monodromy operators. This construction was used
in \cite{TV3} to obtain the explicit formulas for the off-shell
Bethe vectors on the tensor product of evaluation modules of
$\Uqgln$. Two different recurrence relations for the off-shell
Bethe vectors on the evaluation $\Uqgln$-modules were obtained in
\cite{TV3}. Iteration of  these two relations allows to obtain
different explicit formulas for the off-shell Bethe vectors  (see
examples \r{Wt555x} and \r{Wt555y} below). Existence of two types
of the recurrence relations in the nested Bethe ansatz is a
consequence of two different ways of embedding $\Uqgl{N-1}$ into
$\Uqgln$, when they are realized in terms of $\LL$-operators. The
$\LL$-operator of $\Uqgl{N-1}$ can be placed either into the
top-left or into the down-right corners of the $\Uqgln$
$\LL$-operator.

Due to the applications in the theory of quantized
Knizhnik-Zamolodchikov equations the off-shell Bethe vectors are
called {\em weight functions}. We will use both names for these
objects. We will call a weight function {\em universal} if it is
defined in an arbitrary $\Uqgln$-module generated by arbitrary
singular vector.

An alternative approach  for the  construction of the off-shell
Bethe vectors for an arbitrary quantum affine algebra was developed
in \cite{EKhP}. This approach uses the relation between
$\LL$-operator realization of $\Uqgln$ \cite{RS} and the current
realization of the same algebra \cite{D88}. An isomorphism between
these realizations was observed in the paper \cite{DF}. Two
different type of Borel subalgebras are related to these two
realizations  of $\Uqgln$. It was conjectured in \cite{KPT} that
the projections onto intersection of the different type Borel
subalgebras for the product of Drinfeld currents coincide with
the off-shell Bethe vectors obtained using the construction of
\cite{VT}.

The background for the conjecture from \cite{KPT} was an observation that both quantities
satisfy the same coproduct property \cite{EKhP}. It was proved in \cite{KP-GLN} that calculation
of the projection for the product of currents gives similar nested recurrence relation for the off-shell
Bethe vectors as it was obtained in \cite{TV3} on the
level of tensor product of the evaluation $\Uqgln$-modules. This leads to the conclusion that
the projection method yields the universal off-shell Bethe vectors for an
arbitrary $\Uqgln$-module generated by a singuar weight vectors. Only one type of the
recurrence relation which leads to the formula of the type \r{Wt555x} was considered
in the paper \cite{KP-GLN}. Here we generalize the results of this paper. We prove that
in order to get both types of the recurrence relations one has to use two different
isomorphic current realizations of $\Uqgln$. Origin of these different realizations
lies in two possibilities to introduce the Gauss decompositions of $\LL$-operators,
each corresponds to the different embedding of $\Uqgl{N-1}$ $\LL$-operators into
$\Uqgln$ $\LL$-operator (see \r{GF2}--\r{GE2} and \r{GF1}--\r{GE1} below).

The paper is composed as follows. Section~\ref{TVC} serves as reminder of the $\LL$-operator
realization of $\Uqgln$ and construction of the off-shell Bethe vectors in terms of the
matrix elements of these $\LL$-operators \cite{VT}.
In Section~\ref{section2.2} two different Gauss decomposition are introduced as well
as the corresponding current realizations of $\Uqgln$. Here we introduce the current Borel
subalgebras and describe the projections onto intersections of the standard and current
Borel subalgebras. Section~\ref{sec4} devoted to the calculations of the projections
for the product of the currents. The main result here is the Theorem~\ref{main-th} which
yields the universal weight functions as the sum over the ordered products of the projections
of the simple and composed roots currents.
Using the fact  that the projections
of the composed root currents coincide with the Gauss coordinates of $\LL$-operators we
give in the Section~\ref{sec5}
the explicit expressions for the universal weight vectors
in terms of the matrix elements of $\LL$-operators generalizing formulas of the paper
\cite{TV3}.  The main result of the paper is formulated
in the form of the Theorem~\ref{mainmain} in the Section~\ref{sec6}.
The paper contains two Appendices. One contains the reformulation of
the Serre relations in the form of
the commutation relations between composed currents. These relations are necessary to prove
the Proposition~\ref{main-fact}, which is the most difficult technical result of the paper.
It is formulated in the second Appendix and
describes the ordering of the currents and the negative projections of the currents.

\section{Tarasov-Varchenko construction}\label{TVC}

\subsection{$\Uqgln$ in $\LL$-operator formalism}

Let $\E_{ij}\in{\textrm{End}}(\CC^N)$ be a matrix with the only nonzero entry
equal to $1$
at the intersection of the $i$-th row and $j$-th column.
Let $\RR(u,v)\in{\textrm{End}}(\CC^N\ot\CC^N)\ot \CC[[{v}/{u}]]$,
\begin{equation}\label{UqglN-R}
\begin{split}
\RR(u,v)\ =\ &\ \sum_{1\leq i\leq N}\E_{ii}\ot \E_{ii}\ +\ \frac{u-v}{qu-q^{-1}v}
\sum_{1\leq i<j\leq N}(\E_{ii}\ot \E_{jj}+\E_{jj}\ot \E_{ii})
\\
+\ &\frac{q-q^{-1}}{qu-q^{-1}v}\sum_{1\leq i<j\leq N}
(u \E_{ij}\ot \E_{ji}+ v \E_{ji}\ot \E_{ij})
\end{split}
\end{equation}
be a trigonometric $\RR$-matrix associated with the vector
representation of ${\mathfrak{gl}}_N$. Let $q$ be  a complex parameter not equal to zero or
root of unity.

The algebra $\Uqgln$ (with the zero central charge and the gradation operator
dropped out) is an associative algebra with unit generated by the modes
$\LL^{\pm}_{i,j}[\pm k]$, $k\geq 0$, $1\leq i,j\leq N$ of the
$\LL$-operators\footnote{In this paper we use $\RR$-matrix and $\LL$-operators of the
paper \cite{TV3}. See in \cite{KP-GLN} the discussion on the relation between choices
of $\LL$-operators in the papers \cite{TV3} and \cite{KPT}.}
$\LL^{\pm}(z)=\sum_{k=0}^\infty\sum_{i,j=1}^N \E_{ij}\otimes \LL^{\pm}_{i,j}[\pm k]z^{\mp
k}$, subject to the relations
\begin{equation}\label{L-op-com}
\begin{split}
&\RR(u,v)\cdot (\LL^{\pm}(u)\ot \mathbf{1})\cdot (\mathbf{1}\ot \LL^{\pm}(v))=
(\mathbf{1}\ot \LL^{\pm}(v))\cdot (\LL^{\pm}(u)\ot \mathbf{1})\cdot \RR(u,v),\\
&\RR(u,v)\cdot (\LL^{+}(u)\ot \mathbf{1})\cdot (\mathbf{1}\ot \LL^{-}(v))=
(\mathbf{1}\ot \LL^{-}(v))\cdot (\LL^{+}(u)\ot \mathbf{1})\cdot \RR(u,v),
\end{split}
\end{equation}
\begin{equation}\label{L-op-com1}
\LL^{+}_{j,i}[0]=\LL^{-}_{i,j}[0]=0,\qquad \LL^{+}_{k,k}[0]\LL^{-}_{k,k}[0]=1,
\qquad 1\leq i<j \leq N,
\quad 1\leq k\leq N\ .
\end{equation}
Subalgebras formed by the modes $\LL^\pm[n]$ of the $\LL$-operators $\LL^\pm(t)$
are the standard  Borel subalgebras $U_q(\mathfrak{b}^\pm)\subset\Uqgln$.
These Borel subalgebras are Hopf subalgebras of $\Uqgln$.
The coalgebraic structure of these subalgebras is given by the formulae
\begin{equation}\label{coprL}
\Delta \sk{\LL^{\pm}_{i,j}(u)}=\sum_{k=1}^N\ \LL^{\pm}_{k,j}(u)\otimes
\LL^{\pm}_{i,k}(u)\,.
\end{equation}

\subsection{Evaluation homomorphism}

Let $\ele_{a,a}$, $\ele_{a,a+1}$ and $\ele_{a+1,a}$ be
Chevalley generators of the algebra $\Uqqq$ with commutation relations:
\begin{equation*}
\ele_{a,a}\ele_{b,c}\ele^{-1}_{a,a}=q^{\d_{ab}-\d_{ac}}\ele_{b,c}\,,\qquad
[\ele_{a,a+1},\ele_{b+1,b}]=\d_{ab}\frac{\ele_{a,a}\ele^{-1}_{a+1,a+1}-
\ele^{-1}_{a,a}\ele_{a+1,a+1}}{q-q^{-1}}\,,
\end{equation*}
and Serre relations:
\begin{equation}
\begin{split}\label{ser-fin}
\ele^2_{a\pm 1,a}\ele_{a,a\mp 1}-(q+q^{-1})\ele_{a\pm
1,a}\ele_{a,a\mp 1}\ele_{a\pm 1,a}
+\ele_{a,a-1}\ele^2_{a\pm 1,a}&=0\,,\\
\ele^2_{a,a\pm 1}\ele_{a\mp 1,a}-(q+q^{-1})\ele_{a,a\pm 1}\ele_{a\mp
1,a}\ele_{a,a\pm 1}
+\ele_{a\mp 1,a}\ele^2_{a,a\pm 1}&=0\,.
\end{split}
\end{equation}
Evaluation homomorphism of the algebra $\Uqgln$ onto $\Uqqq$ is defined as
\begin{equation}\label{eval2}
\mathcal{E}v_z\sk{{\LL}^+(u)}= {\LL}^+-\, \frac{z}{u}\ {\LL}^-\,,\quad
\mathcal{E}v_z\sk{{\LL}^-(u)}= {\LL}^--\, \frac{u}{z}\ {\LL}^+\,,
\end{equation}
where
\begin{equation*}
\LL^+_{a,b}=\left\{\begin{array}{lc}
(q-q^{-1})\ele_{b,a} \ele_{b,b}&a<b\\
\ele_{a,a}&a=b\\ 0&a>b
\end{array} \right.\,,
\qquad
\LL^-_{a,b}=\left\{\begin{array}{lc}
0&a<b\\
\ele^{-1}_{a,a}&a=b\\ (q^{-1}-q)\ele_{a,a}^{-1} \ele_{b,a}&a>b
\end{array} \right.
\end{equation*}
and the composed roots  generators are defined as follows
\begin{equation*}
\ele_{c,a}= \ele_{c,b}\ele_{b,a}-q^{-1}\ele_{b,a}\ele_{c,b}\,,\qquad
\ele_{a,c}= \ele_{a,b}\ele_{b,c}-q\ele_{b,c}\ele_{a,b}\,,\qquad a<b<c\,.
\end{equation*}
These formulas may be proved inductively after substitution of \r{eval2} into
\r{L-op-com}. The coproduct of the $\LL$-operators \r{coprL}
implies the coproduct of the Chevalley generators: $\Delta\ele_{a,a}=\ele_{a,a}\ot\ele_{a,a}$,
$\Delta\ele_{a,a+1}=\ele_{a,a+1}\ot 1+ \ele_{a,a}^{-1}\ele_{a+1,a+1}\ot \ele_{a,a+1}$ and
$\Delta\ele_{a+1,a}=1\ot\ele_{a+1,a}+\ele_{a+1,a}\ot\ele_{a,a}\ele_{a+1,a+1}^{-1}$.

\subsection{Combinatorial formulas for off-shell Bethe vectors}
\label{subsec2.3}

Let us remind the construction of the off-shell Bethe vectors
 \cite{VT}. Let $\LL$-operator\footnote{We omit superscript + in this $\LL$-operator, since
 will consider only positive standard Borel subalgebra ${U}_q({\mathfrak{b}}^+)$
 here and below.} $\LL^{}(z)=\sum_{k=0}^\infty
\sum_{i,j=1}^N \E_{ij}\otimes \LL^{}_{i,j}[ k]z^{-k}$ of the Borel
subalgebra ${U}_q({\mathfrak{b}}^+)$ of $\Uqgln$ satisfies the Yang-Baxter commutation relation with a
$R$-matrix $\RR(u,v)$.  We use the notation $\LL^{(k)}(z)\in \sk{\CC^N}^{\otimes
M}\ot {U}_q({\mathfrak{b}}^+)$ for
$\LL$-operator acting nontrivially on $k$-th
tensor factor in the product $\sk{\CC^{N}}^{\otimes M}$ for $1\leq k\leq M$.
 Consider a series in $M$ variables
\begin{equation}\label{product}
\ct(u_1,\ldots,u_M)=\LL^{(1)}(u_1)\cdots \LL^{(M)}(u_M)
\cdot \ccR^{(M,\ldots,1)}(u_M,\ldots,u_1)
\end{equation}
with coefficients in $\sk{\End(\CC^N)}^{\ot M}\ot {U}_q({\mathfrak{b}}^+)$,
where
\begin{equation}\label{Rproduct}
\ccR^{(M,\ldots,1)}(u_M,\ldots,u_1)=
{\prodl_{M \geq j > 1}}\ {\prodl_{j > i \geq 1}}\,\,{\RR}^{(ji)}(u_j,u_i)\,.
\end{equation}
In the ordered product of $R$-matrices \r{Rproduct} the factor ${\RR}^{(ji)}$ is to
the left of the factor ${\RR}^{(ml)}$ if $j>m$, or $j=m$ and $i>l$.
Consider the set of  variables
\begin{equation}\label{set111}
\bar{t}_{\segg{\bar{n}}} = \left\{
t^{1}_{1},\ldots,t^{1}_{n_{1}}; t^{2}_{1},\ldots,
t^{2}_{n_{2}}; \ldots\ldots; t^{N-2}_{1},\ldots, t^{N-2}_{n_{N-2}};
t^{N-1}_{1},\ldots,t^{N-1}_{n_{N-1}}\right\}\,.
\end{equation}
Following \cite{VT}, let
\begin{equation}\label{Vit-el}
\begin{split}
&\bbb(\bar t_{\segg{\bar n}})
=\prod_{1\leq a<b\leq N-1}\ \  \prod_{1\leq j\leq n_b}\  \ \prod_{1\leq i\leq n_a}\
\frac{qt^b_j-q^{-1}t^a_i}{t^b_j-t^a_i}\times\\
&\times({\rm tr}_{(\CC^N)^{\ot|{\bar n}|}}\ot {\rm id})\ \big(\ct(t^1_1,\ldots,t^1_{n_1};\ \ldots\ ;
t^{N-1}_1,\ldots,t^{N-1}_{n_{N-1}})\E_{21}^{\ot n_1}\ot\cdots\ot \E_{N\,N-1}^{\ot
n_{N-1}}\ot 1\big),
\end{split}
\end{equation}
where $|\bar n|=n_{1}+\cdots+n_{N-1}$. The element $\ct(\bar t_{\bar n})$
in \r{Vit-el} is given by \r{product}
with identification: $M=|\bar n|$, for $a=1,\ldots,N-1$ and $n_1+\cdots+n_{a-1}<i\leq
n_1+\cdots+n_{a}$, $u_i=t^a_{i-n_1-\cdots-n_{a-1}}$. The coefficients of
$\bbb(\bar t_{\bar n})$ are elements of the Borel subalgebra
$U_q({\mathfrak{b}}^+)$.

For the collection of positive integers $\bar n=\{n_1,\ldots,n_{N-1}\}$ we consider
a direct product of the symmetric groups:
$S_{\bar n} =S_{n_1}\times \cdots
\times S_{n_{N-1}}$.
For any  function $G(\bar t_{\segg{\bar n}})$
we denote by
\begin{equation}\label{symmet}
\Sym_{\ \bar t_{\segg{\bar n}}} \ G(\bar t_{\segg{\bar n}})= \sum_{\si\in
S_{\bar n}}\ G(^\si \bar t_{\segg{\bar n}})
\end{equation}
a symmetrization over groups of variables $\{t^a_1,\ldots,t^a_{n_a}\}$ of the type
$a$, where
\begin{equation}\label{sigmat}
^\si \bar t_{\segg{\bar n}} =\{t^{1}_{\si^{1}(1)},\ldots,
t^{1}_{\si^{1}(n_1)};
\ldots;t^{N-1}_{\si^{N-1}(1)},\ldots,t^{N-1}_{\si^{N-1}(n_{N-1})}\}.
\end{equation}
Let
\begin{equation*}
\qsym(\bar t_{\segg{\bar n}})=\prod_{a=1}^{N-1}\prod_{1\leq \ell<\ell'\leq n_a}
\frac{q^{-1}t^a_\ell-qt^a_{\ell'}}{t^a_\ell-t^a_{\ell'}}
\end{equation*}
be a function of the formal variables $t^a_\ell$.
Following \cite{TV3}, $q$-symmetrization of arbitrary function $G(\bar t_{\segg{\bar n}})$
 means
\begin{equation}\label{qs1r}
\qSym_{\ \bar t_{\segg{\bar n}}} \ G(\bar t_{\segg{\bar n}})=
 \Sym_{\ \bar t_{\segg{\bar n}}} \sk{ \qsym(\bar t_{\segg{\bar n}})\cdot
 G(\bar t_{\segg{\bar n}})}\,.
\end{equation}

We call  vector $v$ \emph{a weight singular vector} if it is annihilated
by any positive mode  $\LL_{i,j}[n]$, $i>j$, $n\geq 0$ of the matrix elements
of the $\LL$-operator $\LL^+(z)$  and is
an eigenvector of the diagonal matrix entries
$\LL^+_{i,i}(z)$
\begin{equation}\label{hwv}
\LL^{+}_{i,j}(z)\ v=0\,,\quad i>j\,,\quad \LL^{+}_{i,i}(z)\ v=\lambda_i(z)\ v\,,\quad
i=1,\ldots,N\,.
\end{equation}
For any $\Uqgln$-module $V$ with a
 singular vector $v$ denote
\begin{equation}\label{Vel-vs-pr}
{\mathbb{B}}_V(\bar t_{\segg{\bar n}})=\bbb(\bar t_{\segg{\bar n}})v.
\end{equation}
Vector valued function ${\mathbb{B}}_{V}(\bar t_{\segg{\bar n}})$
was called in  \cite{VT,TV3} {\it universal off-shell Bethe vector}.

Let $M_\Lambda(z)$ be an evaluation module generated by a singular vector $v$
such that $\ele_{a,a}v=q^{\Lambda_a}v$. From analysis of the relations of the
hierarchical Bethe ansatz \cite{KR83} the authors of the paper \cite{TV3} obtained two
recurrence relations for the off-shell Bethe vectors
${\mathbb{B}}_{M_\Lambda(z)}(\bar t_{\segg{\bar n}})$. Iterating these relations
many equivalent formulas for these objects can be found in terms of the $\Uqqq$
generators $\ele_{a,b}$. Two extreme cases were presented in \cite{TV3}.  First,
the off-shell Bethe vector can be written as
\begin{equation}\label{Wt555x}
\begin{split}
&\ds {\mathbb{B}}_{M_\Lambda(z)}(\bar t_{\segg{\bar n}})=
(q-q^{-1})^{\sum_{a=1}^{N-1}n_a}
\sum_{\admis{s}}\left(\!\!\sk{\
\prod_{N-1\geq b>a\geq 1}^{\longleftarrow}
\frac{q^{s^b_{a-1}(s^b_{a-1}-s^b_{a})}}{[s^b_a-s^b_{a-1}]_q!}
\ \check\ele_{b+1,a}^{s^b_a-s^b_{a-1}}}v\ \right.\\
&\quad\times\ds\left.
\qSym_{\ \bar t_{\segg{\bar n}}}\sk{
\prod_{b=2}^{N-1}\prod_{a=1}^{b-1}\prod_{\ell=1}^{s^b_a}
\frac{q^{\weig_{a+1}}t^{a}_{\ell+\ppp^{b}_{a}}   -q^{-\weig_{a+1}}z }
{t^{a+1}_{\ell+\ppp^{b}_{a+1}}-t^{a}_{\ell+\ppp^{b}_{a}}}
\prod_{\ell'=1}^{\ell+\ppp^b_{a+1}-1}
\frac{qt^{a+1}_{\ell'}-q^{-1}t^{a}_{\ell+\ppp^{b}_{a}}}
{t^{a+1}_{\ell'}-t^{a}_{\ell+\ppp^{b}_{a}}}
}\!\!\right),
\end{split}
\end{equation}
where $\check\ele_{b+1,a}=\ele_{b+1,a}\ele_{b+1,b+1}$ and
sum is taken over all possible collections of nonnegative integers $\admis{\ss}=\{\ss_i^j\}$  such that
\begin{equation}\label{admis-h}
0=\ss^a_0\leq \ss^a_1\leq\cdots\leq \ss^a_a\,,\quad n_a=\sum_{b=a}^{N-1}\ss_a^b\,,\quad a=1,\ldots, N-1
\end{equation}
and $\ppp^b_a=\ss^a_a+\ss^{a+1}_a+\cdots+\ss^{b-1}_a$. Second, the same off-shell Bethe vector
${\mathbb{B}}_{M_\Lambda(z)}(\bar t_{\segg{\bar n}})$ has a different presentation
\begin{equation}\label{Wt555y}
\begin{split}
\ds &{\mathbb{B}}_{M_\Lambda(z)}(\bar t_{\segg{\bar n}})=\ds (q-q^{-1})^{\sum_{a=1}^{N-1}n_a}
\sum_{\admis{\ms}}\left(\!\! \left(\  \prod_{1\leq b\leq a\leq N-1}^{\longrightarrow}
\frac{q^{\ms^b_{a+1}(\ms^b_a-\ms^b_{a+1})}}{[\ms^b_a-\ms^b_{a+1}]_q!}
\check\ele_{a+1,b}^{\ms^b_a-\ms^b_{a+1}} \right)v \right.\\
&\left.\times\ \qSym _{\ \bar t_{\segg{\bar n}}}
\sk{\prod_{a=2}^{N-1}\prod_{b=1}^{a-1}\prod_{\ell=0}^{\ms^b_a-1}
\frac{q^{\weig_a}t^{a}_{\bfm^{b}_{a}-\ell} - q^{-\weig_a} z  }
{t^{a-1}_{\bfm^{b}_{a-1}-\ell}- t^{a}_{\bfm^{b}_{a}-\ell}}
\ \ \prod_{\ell'=\bfm^b_{a-1}-\ell+1}^{n_{a-1}}
\frac{qt^{a}_{\bfm^{b}_{a}-\ell} -q^{-1}t^{a-1}_{\ell'}}
{t^{a}_{\bfm^{b}_{a}-\ell}-t^{a-1}_{\ell'}}}
\!\!\right)\!,
\end{split}
\end{equation}
where $\admis{\ms}=\{\ms_i^j\}$ is a collection nonnegative integers such that
\begin{equation}\label{admis-m}
\ms^a_a\geq \ms^a_{a+1}\geq\cdots\geq \ms^a_{N-1}\geq \ms^a_N=0\,,
\quad n_a=\sum_{b=1}^{a}\ms_a^b\,,\quad a=1,\ldots, N-1
\end{equation}
and $\bfm^{b}_{a}=\ms^1_a+\ms^2_a+\cdots+\ms^b_a$. Ordering product of the non-commutative
entries in  \r{Wt555x} is the same as in \r{Rproduct} and the ordering in \r{Wt555y} is inverse.

More general formula
for the universal off-shell Bethe vectors was obtained in \cite{KP-GLN} using current realization
of the quantum affine algebra $\Uqgln$ and method of projections introduced in
\cite{ER} and developed in  \cite{EKhP}.
Formula \r{Wt555x}  was obtained in the latter paper after specialization to
the evaluation modules. The goal of the present paper is to describe the recurrence relations for the
universal off-shell Bethe vectors or universal weight functions in terms of the
modes of the $\Uqgln$ currents. To do this we have to introduce to different current realizations
of the algebra $\Uqgln$.

\section{Different type Borel subalgebras}
\label{section2.2}

\subsection{Two Gauss decompositions of $\LL$-operators}

The relation between $\LL$-operator realization of $\Uqgln$ and its current realization \cite{D88}
is known since the work \cite{DF}. The main distinction in these two realizations of the same
algebra lies in the different choice of the Borel subalgebras and the corresponding
coalgebraic structures.
To build an isomorphism between two realizations one has to consider
the Gauss decomposition \cite{DF} of the $\LL$-operators and identify linear combinations of certain
Gauss coordinates with the total currents of $\Uqgln$ corresponding to the simple roots of
$\mathfrak{gl}_N$.  In order to construct the universal off-shell Bethe vectors in terms of the
modes of the currents one has to consider the ordered product of the simple roots currents
and calculate the projection of this product onto intersection of the current and standard
Borel subalgebras in $\Uqgln$. Current Borel subalgebras will be introduced in the Subsection~\ref{cBs}.

For the $\LL$-operators fixed by the relations \r{L-op-com} and \r{L-op-com1} we have
two possibilities to introduce the Gauss coordinates $\FF^{\pm}_{b,a}(t)$,
$\EE^{\pm}_{b,a}(t)$, $b>a$ and
$k^\pm_{c}(t)$:
\begin{align}\label{GF2}
\LL^{\pm}_{a,b}(t)&=\FF^{\pm}_{b,a}(t)k^+_{a}(t)+\sum_{1\leq m< a}
\FF^{\pm}_{b,m}(t)k^\pm_{m}(t)\EE^{\pm}_{m,a}(t),\qquad a<b,\\
\label{GK2}
\LL^{\pm}_{a,a}(t)&=k^\pm_{a}(t) +\sum_{1\leq m< a} \FF^{\pm}_{a,m}(t)k^\pm_{m}(t)
\EE^{\pm}_{m,a}(t),\\
\label{GE2}
\LL^{\pm}_{a,b}(t)&=k^\pm_{b}(t)\EE^{\pm}_{b,a}(t)+\sum_{1\leq m<b}
\FF^{\pm}_{b,m}(t)k^\pm_{m}(t)\EE^{\pm}_{m,a}(t),\qquad   a>b\,,
\end{align}
or $\tFF^{\pm}_{b,a}(t)$,
$\tEE^{\pm}_{b,a}(t)$, $b>a$ and
$\tk^\pm_{c}(t)$:
\begin{align}\label{GF1}
\LL^{\pm}_{a,b}(t)&=\tFF^{\pm}_{b,a}(t)\tk^+_{b}(t)+\sum_{b<m\leq N}
\tFF^{\pm}_{m,a}(t)\tk^+_{m}(t)\tEE^{\pm}_{b,m}(t),\qquad a<b,\\
\label{GK1}
\LL^{\pm}_{b,b}(t)&=\tk^\pm_{b}(t) +\sum_{b<m\leq N} \tFF^{\pm}_{m,b}(t)\tk^\pm_{m}(t)
\tEE^{\pm}_{b,m}(t),\\
\label{GE1}
\LL^{\pm}_{a,b}(t)&=\tk^\pm_{a}(t)\tEE^{\pm}_{b,a}(t)+\sum_{a<m\leq N}
\tFF^{\pm}_{m,a}(t)\tk^\pm_{m}(t)\tEE^{\pm}_{b,m}(t),\qquad   a>b\,.
\end{align}

Formulas \r{GF2}--\r{GE1} can be inverted to express the Gauss coordinates
 in terms of the  matrix elements
of $\LL$-operators  as well as to express one set of the Gauss coordinates through another.
This is possible because of the relations:
\begin{equation}\label{car-zmr}
L^\pm_{a,a}[0]=k^\pm_a[0]=\tk^\pm_a[0],\qquad k^+_a[0]\,k^-_a[0]=1
\end{equation}
and  the mode decomposition
 of the Gauss coordinates
\begin{equation}\label{Gau-dec}
\FF_{b,a}^\pm(t)=\sum_{n\geq 0\atop n<0} \FF_{b,a}^\pm[n]t^{-n},\qquad
\EE_{a,b}^\pm(t)=\sum_{n> 0\atop n\leq0} \EE_{a,b}^\pm[n]t^{-n},\quad a<b
\end{equation}
which follows from the mode decomposition of the $\LL$-operators \r{L-op-com1}.
The same rules of the mode decompositions is valid for another set of Gauss coordinates
$\tFF_{b,a}^\pm(t)$ and $\tEE_{a,b}^\pm(t)$.

Let $T$ be a $(n+m)\times(n+m)$ matrix with noncommutative entries presented in the form
\begin{equation*}
 T=\begin{pmatrix}
   A & B \\
   C & D
 \end{pmatrix},
\end{equation*}
where $A$, $B$, $C$ and $D$ are $n\times n$, $n\times m$, $m\times n$ and $m\times m$
matrices respectively.
Suppose that $A$ and $D$ are invertible matrices and introduce the following
notations~\cite{Kl}:
\begin{align}
 \begin{vmatrix}
   \boxed{A} & B \\
   C & D
 \end{vmatrix}&:=A-BD^{-1}C\,, &
 \begin{vmatrix}
   A & B \\
   C & \boxed{D}\,
 \end{vmatrix}&:=\tilde{D}\,, \nonumber
\end{align}
where $\tilde D$ is an $m\times m$ matrix with the entries $\tilde{D}_{ij} =
D_{ij} - \sum\limits_{k,l=1}^n B_{lj} (A^{-1})_{kl}
C_{ik}$. Then the Gauss coordinates can be expressed through the $L$-operator entries as follows
\begin{align}
 k^{\pm}_a(u)&=
 \begin{vmatrix}
   L^{\pm}_{11}(u) & \ldots & L^{\pm}_{1,a-1}(u) & L^{\pm}_{1a}(u) \\
   \vdots &  &\vdots  & \vdots \\
   L^{\pm}_{a-1,1}(u) & \ldots & L^{\pm}_{a-1,a-1}(u) & L^{\pm}_{a-1,a}(u) \\
   L^{\pm}_{a1}(u) & \ldots & L^{\pm}_{a,a-1}(u) & \boxed{L^{\pm}_{aa}(u)}
 \end{vmatrix}, & & \nn \\[2mm]
 \EE^{\pm}_{ab}(u)&=k^{\pm}_a(u)^{-1}
 \begin{vmatrix}
   L^{\pm}_{11}(u) & \ldots & L^{\pm}_{1,a-1}(u) & L^{\pm}_{1a}(u) \\
   \vdots &  &\vdots  & \vdots \\
   L^{\pm}_{a-1,1}(u) & \ldots & L^{\pm}_{a-1,a-1}(u) & L^{\pm}_{a-1,a}(u) \\
   L^{\pm}_{b1}(u) & \ldots & L^{\pm}_{b,a-1}(u) & \boxed{L^{\pm}_{ba}(u)}
 \end{vmatrix}
 , &  &a<b, \nn \\[2mm]
 \FF^{\pm}_{ab}(u)&= \begin{vmatrix}
   L^{\pm}_{11}(u) & \ldots & L^{\pm}_{1,b-1}(u) & L^{\pm}_{1a}(u) \\
   \vdots &  &\vdots  & \vdots \\
   L^{\pm}_{b-1,1}(u) & \ldots & L^{\pm}_{b-1,b-1}(u) & L^{\pm}_{b-1,a}(u) \\
   L^{\pm}_{b1}(u) & \ldots & L^{\pm}_{b,b-1}(u) & \boxed{L^{\pm}_{ba}(u)}
 \end{vmatrix} k^{\pm}_b(u)^{-1}
, &  &a>b.\nn
\end{align}
and
\begin{align}
 \tk^{\pm}_a(u)&=
 \begin{vmatrix}
   \boxed{L^{\pm}_{aa}(u)} & L^{\pm}_{a,a+1}(u) & \ldots & L^{\pm}_{aN}(u) \\
   L^{\pm}_{a+1,a}(u) & L^{\pm}_{a+1,a+1}(u)     & \ldots & L^{\pm}_{a+1,N}(u) \\
   \vdots & \vdots     &  & \vdots \\
   L^{\pm}_{Na}(u) & L^{\pm}_{N,a+1}(u)     & \ldots & L^{\pm}_{NN}(u)
 \end{vmatrix}, & & \nn \\[2mm]
 \tEE_{ab}(u)&=\tk^{\pm}_b(u)^{-1}
 \begin{vmatrix}
   \boxed{L^{\pm}_{ba}(u)} & L^{\pm}_{b,b+1}(u) & \ldots & L^{\pm}_{bN}(u) \\
   L^{\pm}_{b+1,a}(u) & L^{\pm}_{b+1,b+1}(u)     & \ldots & L^{\pm}_{b+1,N}(u) \\
   \vdots & \vdots     &  & \vdots \\
   L^{\pm}_{Na}(u) & L^{\pm}_{N,b+1}(u)     & \ldots & L^{\pm}_{NN}(u)
 \end{vmatrix}, &  &a<b, \nn\\[2mm]
 \tFF^{\pm}_{ab}(u)&=
 \begin{vmatrix}
   \boxed{L^{\pm}_{ba}(u)} & L^{\pm}_{b,a+1}(u) & \ldots & L^{\pm}_{bN}(u) \\
   L^{\pm}_{a+1,a}(u) & L^{\pm}_{a+1,a+1}(u)     & \ldots & L^{\pm}_{a+1,N}(u) \\
   \vdots & \vdots     &  & \vdots \\
   L^{\pm}_{Na}(u) & L^{\pm}_{N,a+1}(u)     & \ldots & L^{\pm}_{NN}(u)
 \end{vmatrix}\tk^{\pm}_a(u)^{-1}, &  &a>b\nn.
\end{align}

Gauss decomposition \r{GF1}--\r{GE1} was used\footnote{Gauss coordinates
$\tFF_{b,a}^\pm(t)$, $\tEE_{a,b}^\pm(t)$, $a<b$ and $\tk_{a}^\pm(t)$ was
denoted in \cite{KP-GLN} as
$\tilde\FF_{b,a}^\pm(t)$, $\tilde\EE_{a,b}^\pm(t)$ and $k_{a}^\pm(t)$ (cf. Sect. 6 in \cite{KP-GLN}).}
in \cite{KP-GLN} in order to obtain
a recurrence relation for the universal weight function \r{uwf} and to prove the conjecture
of \cite{KPT} that the constructions of the off-shell Bethe vectors using $\LL$-operator
approach \cite{VT} and method of projection \cite{EKhP} coincide for arbitrary
$\Uqgln$-module generated by arbitrary singular vectors.
Here we will use another Gauss decomposition \r{GF2}--\r{GE2} in order to get an
alternative recurrence relation for the universal weight function \r{uwf1} which lead
to the formula \r{Wt555y} for the off-shell Bethe vector.

\subsection{The current realization of $\Uqgln$}

Using arguments of the paper \cite{DF} we may obtain for the linear combinations
of the Gauss coordinates
\begin{equation}\label{DF-iso}
F_i(t)=\FF^{+}_{i+1,i}(t)-\FF^{-}_{i+1,i}(t)\,,\quad
E_i(t)=\EE^{+}_{i,i+1}(t)-\EE^{-}_{i,i+1}(t)
\end{equation}
and `diagonal' Gauss coordinates $k^\pm_i(t)$
 the commutation relations of the quantum affine algebra
$\Uqgln$ with the zero
central charge and the gradation operator dropped out.  In terms of the total currents
$F_i(t)$, $E_i(t)$ and the Cartan currents $k^\pm_i(t)$ these commutation relations are
\begin{equation}
\begin{array}{c}\label{gln-com}
(q^{}z-q^{-1}w)E_{i}(z)E_{i}(w)= E_{i}(w)E_{i}(z)(q^{-1}z-q^{}w)\, ,  \\[1em]
(q^{-1}z-qw)E_{i}(z)E_{i+1}(w)= E_{i+1}(w)E_{i}(z)(z-w)\, , \\[1em]
k_i^\pm(z)E_i(w)\left(k_i^\pm(z)\right)^{-1}= \displaystyle\frac{z-w}{q^{-1}z-q^{}w}E_i(w)\, , \\[1em]
k_{i+1}^\pm(z)E_i(w)\left(k_{i+1}^\pm(z)\right)^{-1}=
\displaystyle\frac{z-w}{q^{}z-q^{-1}w}E_i(w)\, , \\[1em]
k_i^\pm(z)E_j(w)\left(k_i^\pm(z)\right)^{-1}=E_j(w), \qquad {\rm if}\quad i\not=j,j+1\, ,  \\[1em]
(q^{-1}z-q^{}w)F_{i}(z)F_{i}(w)= F_{i}(w)F_{i}(z)(q^{}z-q^{-1}w)\, , \\[1em]
(z-w)F_{i}(z)F_{i+1}(w)= F_{i+1}(w)F_{i}(z)(q^{-1}z-qw)\, , \\[1em]
k_i^\pm(z)F_i(w)\left(k_i^\pm(z)\right)^{-1}= \displaystyle\frac{q^{-1}z-qw}{z-w}F_i(w)\, , \\[1em]
k_{i+1}^\pm(z)F_i(w)\left(k_{i+1}^\pm(z)\right)^{-1}=
\displaystyle\frac{q^{}z-q^{-1}w}{z-w}F_i(w)\, , \\[1em]
k_i^\pm(z)F_j(w)\left(k_i^\pm(z)\right)^{-1}=F_j(w), \qquad {\rm if}\quad i\not=j,j+1\, , \\[1em]
[E_{i}(z),F_{j}(w)]= \delta_{{i},{j}}\ \delta(z/w)\
(q-q^{-1})\left(k^-_{i+1}(z)/k^-_{i}(z)-k^+_{i+1}(w)/k^+_{i}(w)\right)\, , \\[1em]
\end{array}
\end{equation}
and the Serre relations
\begin{equation}
\begin{split}
{\rm Sym}_{z_1,z_{2}}
(E_{i}(z_1)E_{i}(z_2)E_{i\pm 1}(w)
&-(q+q^{-1})E_{i}(z_1)E_{i\pm 1}(w)E_{i}(z_2)+\\
&+E_{i\pm 1}(w)E_{i}(z_1)E_{i}(z_2))=0\, ,\\
\label{serre}
{\rm Sym}_{z_1,z_{2}}
(F_{i}(z_1)F_{i}(z_2)F_{i\pm 1}(w)
&-(q+q^{-1})F_{i}(z_1)F_{i\pm 1}(w)F_{i}(z_2)+\\
&+F_{i\pm 1}(w)F_{i}(z_1)F_{i}(z_2))=0\, .
\end{split}
\end{equation}

In order to obtain the commutation relations \r{gln-com} we substitute the
decomposition \r{GF2}-\r{GE2} into commutation relations \r{L-op-com}.
If we substitute into these commutation relations the decomposition
\r{GF1}-\r{GE1}  we obtain
instead of \r{gln-com} slightly different commutation relations for the
total currents
\begin{equation}\label{DF2}
\tiF_i(t)=\tFF^{+}_{i+1,i}(t)-\tFF^{-}_{i+1,i}(t)\,,\quad
\tiE_i(t)=\tEE^{+}_{i,i+1}(t)-\tEE^{-}_{i,i+1}(t)
\end{equation}
and $\tk^\pm_i(t)$ used in the papers \cite{KPT,KP-GLN}:
\begin{gather}
(q^{-1}z-q^{}w)\tiE_{i}(z)\tiE_{i}(w)= \tiE_{i}(w)\tiE_{i}(z)(q^{}z-q^{-1}w)\, ,  \notag \\
(z-w)\tiE_{i}(z)\tiE_{i+1}(w)= \tiE_{i+1}(w)\tiE_{i}(z)(q^{-1}z-qw)\, ,  \notag \\
\tk_i^\pm(z)\tiE_i(w)\left(\tk_i^\pm(z)\right)^{-1}=
\frac{z-w}{q^{-1}z-q^{}w}\tiE_i(w)\, , \notag \\
\tk_{i+1}^\pm(z)\tiE_i(w)\left(\tk_{i+1}^\pm(z)\right)^{-1}=
\frac{z-w}{q^{}z-q^{-1}w}\tiE_i(w)\, , \notag \\
\tk_i^\pm(z)\tiE_j(w)\left(\tk_i^\pm(z)\right)^{-1}=\tiE_j(w),
\qquad {\rm if}\quad i\not=j,j+1\, , \notag \\
 (q^{}z-q^{-1}w)\tiF_{i}(z)\tiF_{i}(w)= \tiF_{i}(w)\tiF_{i}(z)(q^{-1}z-q^{}w)\,
 ,  \label{gln-com1}\\
(q^{-1}z-qw)\tiF_{i}(z)\tiF_{i+1}(w)= \tiF_{i+1}(w)\tiF_{i}(z)(z-w)\, , \notag \\
\tk_i^\pm(z)\tiF_i(w)\left(\tk_i^\pm(z)\right)^{-1}=
\frac{q^{-1}z-qw}{z-w}\tiF_i(w)\, , \notag \\
\tk_{i+1}^\pm(z)\tiF_i(w)\left(\tk_{i+1}^\pm(z)\right)^{-1}=
\frac{q^{}z-q^{-1}w}{z-w}\tiF_i(w)\, , \notag \\
\tk_i^\pm(z)\tiF_j(w)\left(\tk_i^\pm(z)\right)^{-1}=\tiF_j(w),
 \qquad {\rm if}\quad i\not=j,j+1\, , \notag \\
[\tiE_{i}(z),\tiF_{j}(w)]= \delta_{{i},{j}}\ \delta(z/w)\
(q-q^{-1})\left( \tk^+_{i}(z)/\tk^+_{i+1}(z)-\tk^-_{i}(w)/\tk^-_{i+1}(w)\right)\, \notag
\end{gather}
and the same Serre relations \r{serre} with currents $\ee_{i}(z)$ and $\ff_{i}(z)$
replaced by $\tiE_{i}(z)$ and $\tiF_{i}(z)$.

Formulae \eqref{gln-com} and \eqref{gln-com1} should be considered as
formal series identities describing the infinite set of the relations between modes of the
currents. The symbol $\delta(z)$ entering these relations is a formal series $\sum_{n\in\ZZ} z^n$.
These commutation relations describes two isomorphic
current realizations of the algebra $\Uqgln$. Isomorphism follows from the result of \cite{DF}.

{}For any series $\acc(t)=\sum_{m\in\ZZ}\acc[m]t^{-m}$ we denote
$\acc(t)^{(+)}=\sum_{m>0} \acc[m]\, t^{-m}\,,
$ and $\acc(t)^{(-)}=-\sum_{ m\leq 0} \acc[m]\, t^{-m}\,.\
$
The initial conditions \rf{L-op-com1} imply the relations
\begin{equation}\label{DFinverse}
\FF^{\pm}_{i+1,i}(z)=z\left(z^{-1}F_i(z)\right)^{(\pm)},\qquad
\EE^{\pm}_{i,i+1}(z)=E_i(z)^{(\pm)}
\end{equation}
and the same for the relations between Gauss coordinates
$\tEE^{\pm}_{i,i+1}(t)$, $\tFF^{\pm}_{i+1,i}(t)$ and the currents
$\tiE_i(t)$, $\tiF_i(t)$.

\subsection{Borel subalgebras and projections to their intersections}
\label{cBs}

We consider two types of Borel subalgebras of the
algebra $\Uqgln$. Borel subalgebras $U_q(\mathfrak{b}^\pm)\subset \Uqgln $
are generated by the modes of the $\LL$-operators $\LL^{(\pm)}(z)$ respectively.
For the generators in these subalgebras we can use either modes of the Gauss coordinates
\r{GF2}--\r{GE2},
$\EE^\pm_{i,i+1}(t)$, $\FF^\pm_{i+1,i}(t)$, $k^\pm_{j}(t)$ or the modes
of the Gauss coordinates \r{GF1}--\r{GE1}
$\tEE^\pm_{i,i+1}(t)$, $\tFF^\pm_{i+1,i}(t)$, $\tk^\pm_{j}(t)$.

Another types of Borel subalgebras are related to the current realizations of
$\Uqgln$ given in the previous subsection.
We consider first the  current Borel subalgebras generated by the modes
of the currents $E_i(t)$, $F_i(t)$, $k^\pm_j(t)$.

 The Borel subalgebra $U_F\subset \Uqgln$ is generated by
modes of the currents
$F_i[n]$, $k^+_j[m]$, $i=1,\ldots,N-1$, $j=1,\ldots,N$, $n\in\ZZ$
and $m\geq0$. The Borel subalgebra
$U_E\subset \Uqgln$ is generated by modes of the currents
$E_i[n]$, $k^-_j[-m]$, $i=1,\ldots,N-1$, $j=1,\ldots,N$, $n\in\ZZ$ and
$m\geq0$. We will consider also a subalgebra $U'_F\subset U_F$, generated
by the elements
$F_i[n]$, $k^+_j[m]$, $i=1,\ldots,N-1$, $j=1,\ldots,N$, $n\in\ZZ$
and $m>0$, and a subalgebra $U'_E\subset U_E$ generated by
the elements
$E_i[n]$, $k^-_j[-m]$, $i=1,\ldots,N-1$, $j=1,\ldots,N$, $n\in\ZZ$
and $m>0$.
Further, we will be interested in the intersections,
\begin{equation}\label{Intergl}
U_f^-=U'_F\cap U_q(\mathfrak{b}^-)\,,\qquad
U_F^+=U_F\cap U_q(\mathfrak{b}^+)\,
\end{equation}
and will describe properties of projections to these intersections.
We call $U_F$ and $U_E$ {\em the current Borel subalgebras}.
Let $U_f\subset U_F$ be subalgebra of the current Borel subalgebra generated
by the modes of the currents
$F_i[n]$, $i=1,\ldots,N-1$, $n\in\ZZ$ only. In what follows we will use also
the subalgebra $U^+_f\subset U_f$ defined by the intersection
\begin{equation}\label{Uf+}
U_f^+=U^+_F\cap U_f\,.
\end{equation}

 In \cite{D88} the \emph{current} Hopf structure for the algebra
$\Uqgln$ has been defined,
\begin{equation}\label{gln-copr}
\begin{split}
\Delta^{(D)}\sk{E_i(z)}&=1\ot E_i(z) + E_i(z)\ot k^-_{i+1}(z)\sk{k^-_{i}(z)}^{-1},\\
\Delta^{(D)}\sk{F_i(z)}&=F_i(z)\ot 1 +
k^+_{i+1}(z)\sk{k^+_{i}(z)}^{-1}\ot F_i(z),\\
\Delta^{(D)}\sk{k^\pm_i(z)}&=k^\pm_i(z)\ot k^\pm_{i}(z).
\end{split}
\end{equation}
With respect to the current Hopf structure the current Borel subalgebras are Hopf
subalgebras of $\Uqgln$.
We can check \cite{EKhP,KPT} that  the intersections  $U_f^-$
and $U_F^+$ are subalgebras and coideals with respect to the Drinfeld coproduct
\r{gln-copr}
\begin{equation*}
\Delta^{(D)}(U_F^+)\subset  U_F^+\ot\Uqgln\,,\qquad
\Delta^{(D)}(U_f^-)\subset \Uqgln\ot U_f^- \,,
\end{equation*}
and the multiplication $m$ in $\Uqgln$ induces an isomorphism of vector spaces
$$m: U_f^-\ot U_F^+\to U_F\,.$$
According to the general theory presented in \cite{EKhP} we  define
the projection operators $\Pfp:U_F\subset \Uqgln \to U_F^+$ and
$\Pfm:U_F\subset \Uqgln \to U_f^-$ by the prescriptions
\begin{equation}\label{pgln}
\begin{split}
\Pfp(\Ff_-\ \Ff_+)&=\coun(\Ff_-)\ \Ff_+, \qquad
\Pfm(\Ff_-\ \Ff_+)=\Ff_-\ \coun(\Ff_+),
\\& \text{for any}\qquad \Ff_-\in U_f^-,
\quad \Ff_+\in U_F^+ .
\end{split}
\end{equation}

Denote by  $\overline U_F$ an extension of the algebra $U_F$ formed
by infinite sums of monomials which are ordered products
$a_{i_1}[n_1]\cdots a_{i_k}[n_k]$ with $n_1\leq\cdots\leq n_k$,
where  $a_{i_l}[n_l]$ is either $F_{i_l}[n_l]$ or $k^+_{i_l}[n_l]$.
It was proved in \cite{EKhP} that
\begin{itemize}
\item[(1)] the action of the projections \r{pgln} can be extended to the
 agebra
$\overline U_F$;
\item[(2)] for any $\Ff\in \overline U_F$ with $\Delta^{(D)}(\Ff)=\sum_i \Ff'_i\otimes \Ff''_i$ we have
\begin{equation}\label{pr-prop}
\Ff=\sum_i\Pfm(\Ff'_i)\cdot \Pfp(\Ff''_i)\,.
\end{equation}
\end{itemize}

In \cite{KPT,KP-GLN}, we used the current Borel subalgebras $\hU_F$, $\hU_E$
generated by the modes of the currents $\tiF_i(t)$, $\tiE_i(t)$, $\tk^\pm_j(t)$
in the same way as it was done above for $U_F$, $U_E$. These current Borel subalgebras
are Hopf subalgebras of $\Uqgln$ with a
different Drinfeld coproduct
\begin{equation}\label{gln-cm}
\begin{split}
\hat\Delta^{(D)}\sk{\tiE_i(z)}&= \tiE_i(z)\ot1 + \tk^-_{i}(z)\sk{\tk^-_{i+1}(z)}^{-1}\ot \tiE_i(z)\,,\\
\hat\Delta^{(D)}\sk{\tiF_i(z)}&=1\ot \tiF_i(z) +
\tiF_i(z) \ot \tk^+_{i}(z)\sk{\tk^+_{i+1}(z)}^{-1}\,,\\
\hat\Delta^{(D)}\sk{\tk^\pm_i(z)}&=\tk^\pm_i(z)\ot \tk^\pm_{i}(z)\,.
\end{split}
\end{equation}
The standard Borel subalgebras $U_q({\mathfrak{b}}^{\pm})$ are defined by the modes of
the Gauss coordinates $\tEE^\pm_{i,i+1}(t)$, $\tFF^\pm_{i+1,i}(t)$, $\tk^\pm_{j}(t)$ and their
intersections with the currents Borel subalgebras $\hU_F$ are defined by the same
 formulas  \r{Intergl}.
Using coproduct \r{gln-cm}
one may check that the coalgebraic properties of these intersections are changed to
\begin{equation}\label{hatint}
\hat\Delta^{(D)}(\hU_F^+)\subset \Uqgln\ot \hU_F^+\,,\qquad
\hat\Delta^{(D)}(\hU_f^-)\subset \hU_f^-\ot \Uqgln\,.
\end{equation}
Denote by  $\hat{\overline{U}}_F$ an extension of the algebra $\hat U_F$ formed
by infinite sums of monomials which are ordered products
$a_{i_1}[n_1]\cdots a_{i_k}[n_k]$ with $n_1\leq\cdots\leq n_k$,
where  $a_{i_l}[n_l]$ is either $\hat F_{i_l}[n_l]$ or $\hat k^+_{i_l}[n_l]$.
Projections to the intersections \r{hatint} are defined by the formulas analogous to
\r{pgln} and
can be extended to the algebra  $\hat{\overline{U}}_F$. The
property \r{pr-prop} is  changed to
$\Ff=\sum_i\hPfm(\Ff''_i)\cdot \hPfp(\Ff'_i)$,
where $\hat\Delta^{(D)}(\Ff)=\sum_i \Ff'_i\otimes \Ff''_i$.

\section{Universal weight function and projections}\label{sec4}

We will use the same notations as in \cite{KP-GLN}.

Let $\Pi$ be the set $\{1,\ldots,N-1\}$ of indices of the simple
positive roots of $\mathfrak{gl}_N$. A finite collection
$I=\{i_1,\dots,i_n\}$ with a linear ordering $i_i\ord\cdots\ord i_n$
and a  map $\iota:I\to\Pi$ is called an {\it ordered
$\Pi$-multiset}.
To each $\Pi$-ordered multiset $I=\{i_1,\dots,i_n\}$ we attach an
ordered set of variables $\{t_i|i\in I\}=\{t_{i_1},\dots,t_{i_n}\}$.
Each element $i_k\in I$ and each variable $t_{i_k}$ has its own
`type': $\iota(i_k)\in\Pi$.

Our basic calculations are performed on a level of formal series attached
 to certain ordered multisets. For the save of space we often write some series
as rational homogeneous functions with the following prescription.
Let  $\{t_i\mid i\in I\}=\{t_{i_1},\dots,t_{i_n}\}$ be the ordered set of
variables attached to an ordered set $I=\{i_1\prec
i_2\prec\cdots\prec i_n\}$ and $g(t_i\mid i\in I)$ be a rational
function. Then we associate to $g(t_i\mid i\in I)$ a Loran series which
is the expansion of $g(t_i\mid i\in I)$ in a region
 $|t_{i_1}|\ll|t_{i_2}|\ll\cdots\ll|t_{i_n}|$. If, for instance,
$1\prec 2$, then we associate to a rational function $(t_1-t_2)^{-1}$ a series
$-\sum_{k\geq 0}t_1^k t_2^{-k-1}$. On the contrary, for the ordering $2\prec 1$ we associate to the
same rational function $(t_1-t_2)^{-1}$ a series
$\sum_{k\geq 0}t_2^k t_1^{-k-1}$.

Let $\bar\lll$ and $\bar\rr$ be two collections of nonnegative integers
satisfying a set of inequalities
\begin{equation}\label{set23}
\lll_a\leq\rr_a\,,\quad a=1,\ldots,N-1\,.
\end{equation}
Denote by $\meg{\bar\lll}{\bar\rr}$  a set of segments which contain  positive
integers $\{\lll_a+1,\lll_a+2,\ldots,\rr_a-1,\rr_a\}$
including $\rr_a$ and excluding $\lll_a$. The length of each segment is equal
to $\rr_a-\lll_a$.

For a given set $\meg{\bar\lll}{\bar\rr}$ of segments  we denote by
$\bar t_{\meg{\bar\lll}{\bar\rr}}$  the sets of variables
\begin{equation}\label{set2m}
\bar{t}_{\meg{\bar{\lll}}{\bar\rr}}\
=\{t^{1}_{\lll_{1}+1},\ldots,t^{1}_{\rr_{1}};
t^{2}_{\lll_{2}+1},\ldots,t^{2}_{\rr_{2}};\ldots;
t^{N-1}_{\lll_{N-1}+1},\ldots,t^{N-1}_{\rr_{N-1}} \}\,.
\end{equation}
The number of the variables of the type $a$  is equal to $\rr_a-\lll_a$.
In this notation, the set of variables \r{set111} is $\bar
t_{\segg{\bar n}}\equiv \bar t_{\seg{\bar n}{\bar 0}}$. We will
consider \rf{set2m} as a list of ordered variables, corresponding to two
ordered multisets:
\begin{equation}\label{sord}
I=\{r_{N-1}\prec\cdots \prec  l_{N-1}+1\prec\cdot\ \cdot\prec r_2\prec\cdots\prec l_2+1\prec
    r_1\prec\cdots\prec l_1+1\}
\end{equation}
and
\begin{equation}\label{hsord}
\hat I=\{l_1+1\prec\cdots\prec r_1\prec l_2+1\prec\cdots\prec r_2\prec\cdot\
\cdot\prec l_{N-1}+1\prec\cdots \prec r_{N-1}\}\,.
\end{equation}

{}For any $a=1,\ldots,N-1\,$ we denote the sets of variables corresponding to the  segments
$\meg{\lll_a}{\rr_a}=\{\lll_a+1,\lll_a+2,\ldots,\rr_a\}$
as $\bar{t}^a_{\meg{{\lll_a}}{\rr_a}}\
=\{t^{a}_{\lll_{a}+1},\ldots,t^{a}_{\rr_{a}} \}$.
All the variables in $\bar{t}^a_{\meg{{\lll_a}}{\rr_a}}$ have  the type $a$. For the segments
$\meg{\lll_a}{\rr_a}=\meg{0}{n_a}$ we use the shorten  notations
$\bar{t}_{\meg{\bar 0}{\bar n}}\equiv\bar{t}_{\megg{\bar n}}$ and
$\bar{t}^a_{\meg{0}{n_a}}\equiv\bar{t}^a_{\megg{n_a}}$.

For a collection of variables $\bar t_{\seg{\bar\rr}{\bar\lll}}$ we
consider two types of the ordered products of the currents
\begin{equation}\label{FFFl}
\F(\bar t_{\seg{\bar\rr}{\bar\lll}})=\!\!\prod_{1\leq a\leq N-1}
^{\longrightarrow} \sk{\prod_{\lll_a<\ell\leq\rr_a}^{\longrightarrow} F_{a}(t^a_{\ell})}=
F_{1}(t^1_{\lll_1+1})\cdots F_{1}(t^{1}_{\rr_{1}})\cdots
F_{N-1}(t^{N-1}_{\rr_{N-1}})
\end{equation}
and
\begin{equation}\label{tFFFl}
\hF(\bar t_{\seg{\bar\rr}{\bar\lll}})=\!\!\prod_{N-1\geq a\geq
1}^{\longleftarrow} \sk{\prod_{\rr_a\geq\ell>
\lll_a}^{\longleftarrow} \tiF_{a}(t^a_{\ell})}=
\tiF_{N-1}(t^{N-1}_{\rr_{N-1}})\cdots
\tiF_{1}(t^{1}_{\rr_{1}})\cdots \tiF_{1}(t^1_{\lll_1+1})\,,
\end{equation}
where the series $F_a(t)\equiv F_{a+1,a}(t)\,$ and $\tiF_a(t)\equiv \tiF_{a+1,a}(t)\,$ are
 defined by \r{DF-iso} and \r{DF2}, respectively.
 As particular cases, we have $\,\F(\bar t^a_{\seg{\bar\rr_a}{\bar\lll_a}})=
F_{a}(t^a_{l_a+1}) F_{a}(t^a_{l_a+2})\cdots F_{a}(t^a_{r_a})$ and
$\,\hF(\bar t^a_{\seg{\bar\rr_a}{\bar\lll_a}})=
\tiF_{a}(t^a_{r_a})\cdots \tiF_{a}(t^a_{l_a+2})
\tiF_{a}(t^a_{l_a+1})$.

 Symbols $\mathop{\prod}\limits^{\longleftarrow}_a A_a$ and
$\mathop{\prod}\limits^{\longrightarrow}_a A_a$
 mean  ordered products of
noncommutative entries $A_a$, such that $A_a$ is on the right (resp., on the left)
 from $A_b$ for $b>a$:
\begin{equation*}
\mathop{\prod}\limits^{\longleftarrow}_{j\geq a\geq i} A_a = A_j\,A_{j-1}\,
\cdots\, A_{i+1}\,A_i\,,\quad
\mathop{\prod}\limits^{\longrightarrow}_{i\leq a\leq j} A_a = A_i\,A_{i+1}\,
\cdots\, A_{j-1}\,A_j\,.
\end{equation*}

According to \cite{DKh,E,EKhP,KP} the product of the currents \r{FFFl}
is a formal series over the ratios $t^b_k/t^c_l$ with $b>c$ and $t^a_i/t^a_j$ with $i>j$
 taking values in the completions
$\overline U_F$.
Analogously, the product of the currents \r{tFFFl}
is a formal series over the ratios
 $t^b_k/t^c_l$ with $b<c$ and $t^a_i/t^a_j$ with $i<j$
 taking values in the completion $\hat{\overline{U}}_F$.

It means that these products have
the same analytical structure as the rational functions of the variables
$\bar t_{\seg{\bar\rr}{\bar\lll}}$
defined by the multisets $I$ and $\hat I$, respectively.
The domains of analyticity of the rational functions defined by the multisets
$I$ and $\hat I$ are different. The products \r{FFFl}  and \r{tFFFl} have poles
for some values of the ratios $t^b_k/t^c_l$ and $t^a_i/t^a_j$. The operator valued
coefficient at these poles take values in the completions $\overline U_F$  and
$\hat{\overline{U}}_F$ and can be identified with composed root currents (see
\cite{KP,KP-GLN}).
In what follows we will consider  projections of the product of the currents
\begin{equation}\label{uwf1}
\calW^{N}(\bar{t}_{\segg{\bar{n}}})=
\Pfp\left(F_{1}(t_{1}^{1})\cdots F_{1}(t_{n_{1}}^{1})
  \quad \cdots\quad
F_{N-1}(t_{1}^{N-1})\cdots F_{N-1}(t_{n_{N-1}}^{N-1})
\right)
\end{equation}
and
\begin{equation}\label{uwf}
\ticalW^{N}(\bar{t}_{\segg{\bar{n}}})= \hPfp\left(
\tiF_{N-1}(t_{n_{N-1}}^{N-1})\cdots \tiF_{N-1}(t_{1}^{N-1})
  \quad \cdots\quad
\tiF_{1}(t_{n_1}^1)\cdots \tiF_{1}(t_{1}^1)
\right)\,.
\end{equation}
It was proved in \cite{KP} that these projections can be analytically continued from
there domains of definitions. This allows to compare the universal weight functions defined by these
projections.

It was conjectured in \cite{KPT}  and then proved in \cite{KP-GLN}
that   the universal weight function can be identified
with the projection \r{uwf}.
A method of  computation of this projection
was proposed in \cite{KP} and it was further developed in \cite{KP-GLN}.
In this paper we will calculate the universal weight function
given by the projection \r{uwf1}.

 For any weight singular vector $v$ \r{hwv} let
the related weight functions
\begin{equation}\label{b-v-sym1}
{\bf w}_V^{N}(\bar t_{\segg{\bar n}})=\qsym(\bar t_{\segg{\bar n}})
\ \calW^{N}(\bar{t}_{\segg{\bar{n}}})
\prod_{a=1}^{N-1}\prod_{\ell=1}^{n_{a}}k^+_{a}(t^a_\ell)\ v
\end{equation}
\begin{equation}\label{b-v-sym}
\hat{\bf w}_V^{N}(\bar t_{\segg{\bar n}})=\qsym(\bar t_{\segg{\bar n}})\
\ticalW^{N}(\bar{t}_{\segg{\bar{n}}})\prod_{a=1}^{N-1}\prod_{\ell=1}^{n_{a}}\tk^+_{a+1}(t^{a}_\ell)\ v\,,
\end{equation}
be the functions taking values in the $\Uqgln$-module $V$ generated by a singular vector
$v$. In \cite{KPT,KP-GLN} they were called a modified weight function or
{\em  universal off-shell Bethe vector}.

We will show that analytical continuations of the universal off-shell Bethe vectors defined
by the universal weight functions \r{uwf} and \r{uwf1} coincide:
\begin{equation}\label{thefact}
{\bf w}_V^{N}(\bar t_{\segg{\bar n}})=\hat{\bf w}_V^{N}(\bar t_{\segg{\bar n}})\,.
\end{equation}
This will follow from the fact that
\begin{equation*}
\hat{\bf w}_V^{N}(\bar t_{\segg{\bar n}})={\mathbb{B}}_V(\bar t_{\segg{\bar
n}})
\end{equation*}
for arbitrary $\Uqgln$-module $V$ generated by a singular vector $v$ \cite{KP-GLN}.
In this paper we will prove that
\begin{equation}\label{this-pr}
{\bf w}_V^{N}(\bar t_{\segg{\bar n}})={\mathbb{B}}_V(\bar t_{\segg{\bar n}})\,,
\end{equation}
thus yielding the equality \r{thefact}.

To prove \r{this-pr} we use the same arguments as in  \cite{KP-GLN}. We
calculate the projection \r{uwf1}, rewrite the corresponding universal off-shell Bethe
vector \r{b-v-sym1} in terms of the ordered products of the
matrix elements of $\LL$-operators acting onto singular vector $v$ and show that this  implies
an expression \r{Wt555y} for this vector on the evaluation modules.
This gives a different formula than obtained
in \cite{KP-GLN} for off-shell Bethe vectors.
Using the result of \cite{EKhP} we may check that the universal off-shell Bethe
vector \r{b-v-sym1} defined by the projection satisfies the same comultiplication
properties as the vector ${\mathbb{B}}_V(\bar t_{\segg{\bar n}})$ (see \cite{KPT}). This means
that equality \r{this-pr} is valid for arbitrary tensor product of the
evaluation representations of $\Uqgln$. Then the classical result \cite{CP} implies
that  \r{this-pr} is true for every irreducible finite-dimensional $\Uqgln$-module
$V$ generated by a singular vector $v$.

We conclude this subsection by describing our strategy of calculation of the projections.
We will use the same approach to calculate \r{uwf} that
was used in \cite{KP-GLN} for the calculation of \r{uwf1}.
 In that paper we  separate all  factors $\tiF_a(t^a_i)$ with $a<N-1$
in \r{uwf} and apply to this
product the ordering procedure  based on the
property \r{dec-ff1m}. Here we will separate the factors $F_a(t^a_i)$ with $a>1$ in \r{uwf1}.
In both cases we get under total
projection a symmetrization of a sum of terms
$x_iP^-(y_i)P^+(z_i)$ with rational coefficients. Here
$x_i$ are expressed via modes of $\tiF_{N-1}(t)$, and $y_i$, $z_i$ via modes
of $\tiF_a(t)$ with $a<N-1$ in the case of \r{uwf}. In the case of \r{uwf1},
$x_i$ are expressed via modes of $F_{1}(t)$, and $y_i$, $z_i$ via modes
of $F_a(t)$ with $a>1$.
As well as in \cite{KP-GLN}
we reorder $x_i$ and $P^-(y_i)$ in both cases. At this
stage  composed currents of different types collected in so called strings appear.

\subsection{Symmetrization}

Consider permutation group $S_n$ and its action on the formal series of $n$ variables
defined on the elementary transpositions $\sigma_{i,i+1}$ as follows
\begin{align*}
\pi(\sigma_{i,i+1})G(u_1, \dots,u_i, u_{i+1},\dots, u_n) =
\frac{qu_i-q^{-1}u_{i+1}}{q^{-1}u_i-qu_{i+1}}\ G(u_1, \dots,u_{i+1}, u_i, \dots, u_n).
\end{align*}
The $q$-depending factor in this formula is chosen in such a way that the products
$F_a(u_1)\cdots F_a(u_n)$ and $\tiF_a(u_n)\cdots \tiF_a(u_1)$
are invariant under this action. Summing the action over all the group of permutations
we obtain the operator $\tSym_u=\sum\limits_{\sigma\in S_n}\pi(\sigma)$ acting as follows
\begin{align}\label{qss}
\tSym_u G(u) =\sum\limits_{\si \in S_n}\prod\limits_{\substack{\ell<\ell'\\ \si(\ell)>\si(\ell')}}
\frac{qt_{\si(\ell')}-q^{-1}t_{\si(\ell)}}
{q^{-1}t_{\si(\ell')}-q^{}t_{\si(\ell)}}\  G(^\sigma u).
\end{align}
The product is taken over all pairs $(\ell, \ell')$, such that conditions $\ell
< \ell'$ and $\si(\ell) > \si(\ell')$ are satisfied simultaneously.
The indices of the set of formal variables $u=(u_1,\ldots,u_n)$
are from two ordered multisets: $\{1\prec 2\prec\cdots\prec n\}$ or
$\{n\prec n-1\prec\cdots\prec 1\}$.
The expansions of the rational functions entering the right hand side of \r{qss}
are defined  by these ordered multisets.

We call operator  $\tSym_u$  -- {\em $q$-symmetrization}.
This operation differs from the ones used in the paper \cite{TV3}.
The exact relation between them will be given below.
The operator $\frac1{n!}\tSym_u$ is the group average
with respect to the action $\pi$, which implies the
formula
\begin{equation}\label{sym*}
\tSym_u\tSym_u(\cdot)=n!\tSym_u(\cdot)\,.
\end{equation}

An important property of $q$-symmetrization is the formula
\begin{align}
 \tSym_{(u_1,\ldots,u_n)}
  =\sum_{\sigma\in
    S^{(s)}_n}\pi(\sigma)\;\tSym_{(u_1,\ldots,u_s)}\tSym_{(u_{s+1},\ldots,u_n)}\,,
    \label{sym_div}
\end{align}
where $s\in\seg{0}{n}$ is fixed and the sum is taken over the subset
\begin{align*}
 S^{(s)}_n=\{\sigma\in S_n \mid \sigma(1)<\ldots<\sigma(s),\,\sigma(s+1)<\ldots<\sigma(n)\}\,.
\end{align*}

As in Subsection~\ref{subsec2.3} we denote by $S_{\bar\lll,\bar\rr} =
S_{\lll_{1},\rr_{1}}\times \cdots \times S_{\lll_{N-1},\rr_{N-1}}$
 a direct product of the groups $S_{\lll_{a},\rr_{a}}$ permuting
integer numbers $\lll_{a}+1,\ldots, \rr_{a}$.
The $q$-sym\-me\-tri\-za\-tion over whole set of variables
$\bar t_{\seg{\bar\rr}{\bar\lll}}$ is defined by the formula
\begin{equation}\label{qsr}
\tSym_{\ \bar t_{\seg{\bar\rr}{\bar\lll}}} \ G(\bar
t_{\seg{\bar\rr}{\bar\lll}})= \sum_{\si\in
S_{\bar\lll,\bar\rr}}\prod_{1\leq a\leq N-1}
\prod_{\substack{\ell<\ell'\\ \si^a(\ell)>\si^a(\ell')}}
\frac{q^{}t^a_{\si^a(\ell')}-q^{-1}t^a_{\si^a(\ell)}}
{q^{-1}t^a_{\si^a(\ell')}-q^{}t^a_{\si^a(\ell)}}\ G(^\si \bar
t_{\seg{\bar\rr}{\bar\lll}})\,,
\end{equation}
where  the set $^\si \bar t_{\seg{\bar\rr}{\bar\lll}}$
is defined in the same way  as in \r{sigmat}.

The $q$-symmetrization~\eqref{qsr} is related to the $q$-symmetrization~\r{qs1r}:
\begin{equation}\label{relat}
\qSym_{\ \bar t_{\seg{\bar\rr}{\bar\lll}}} \ G(\bar
t_{\seg{\bar\rr}{\bar\lll}})=
\prod_{a=1}^{N-1}\ \ \prod_{\lll_a<\ell<\ell'\leq\rr_a}\ \ \frac{q^{-1}t^a_\ell-qt^a_{\ell'}}
{t^a_\ell-t^a_{\ell'}}\ \
\tSym_{\ \bar t_{\seg{\bar\rr}{\bar\lll}}} \ G(\bar
t_{\meg{\bar\rr}{\bar\lll}})
\end{equation}

As in \cite{KP-GLN} we  use the $q$-symmetrization \r{qsr} to compute the recurrence
relation for the universal weight function \r{uwf1} and will restore the $q$-symmetrization \r{qs1r}
when express the universal off-shell Bethe vector in terms of matrix elements of $\LL$-operators.

We say that the series $Q(\bar t_{\seg{\bar\rr}{\bar\lll}})$ is $q$-symmetric, if
it is invariant under the action $\pi$ of each group
$S_{\lll_{a},\rr_{a}}$ with respect to the permutations of the variables
$t^a_{\lll_a+1},\ldots,t_{\rr_a}$ for $a=1,\ldots,N-1$. The $q$-symmetrization of $q$-symmetric series
is equivalent to the multiplication by the order of the group $S_{\bar\lll,\bar\rr}$:
\begin{equation}\label{exa3}
\tSym_{\ \bar t_{\seg{\bar\rr}{\bar\lll}}}
Q(\bar t_{\seg{\bar\rr}{\bar\lll}})=\prod_{a=1}^{N-1} (\rr_a-\lll_a)!\
Q(\bar t_{\seg{\bar\rr}{\bar\lll}})\,.
\end{equation}
The $q$-symmetrization $Q(\bar t_{\seg{\bar\rr}{\bar\lll}})=
\tSym_{\bar t_{\seg{\bar\rr}{\bar\lll}}} G(\bar t_{\seg{\bar\rr}{\bar\lll}})$
of any series $ G(\bar t_{\seg{\bar\rr}{\bar\lll}})$ is a $q$-symmetric series,
which follows from \r{sym*}.

Let $G^{\rm sym}(u_1,\ldots,u_n)$ be symmetric function of $n$
 variables $u_k$, i.e.  $G^{\rm sym}(^\si\bar u)=G^{\rm sym}(\bar u)$
for any element $\si$ from the symmetric group $S_n$. Then one can check the
following property of $q$-symmetrization:
\begin{align}\label{po-sim}
\frac{1}{n!}\ \ \tSym_{\ \bar u}\sk{\qsym(\bar u)^{-1}G^{\rm sym}(\bar u)}
&=\frac{1}{[n]_q!}\ \  \tSym_{\ \bar u}\sk{G^{\rm sym}(\bar u)},\\
\qsym(\bar u)&=\prod\limits_{k<k'}
\frac{q^{-1}u_{k}-qu_{k}}{u_k-u_{k'}}\,,\label{betadef}
\end{align}
where $[n]_q=\frac{q^n-q^{-n}}{q-q^{-1}}$ and $[n]_q!=[n]_q[n-1]_q\cdots [2]_q[1]_q$.

\subsection{Normal ordering of the products of currents}

We call any expression $\sum_i f^{(i)}_-\cdot f^{(i)}_+$, where
$f^{(i)}_-\in U^-_f$ and $f^{(i)}_+\in U^+_f$ {\it (normal) ordered}.
Subalgebras $U^-_f$ and  $U^+_f$ are defined by \r{Intergl} and \r{Uf+} respectively.
Using the property of the projections \r{pr-prop}
we can present any product of the currents in a normal ordered form.
Taking into account property~\eqref{sym_div} we obtain the following
\begin{proposition}\label{decomp} We have formal series equalities
\begin{equation}\label{dec-ff1}
\begin{array}{c}
\ds \F(\bar t_{\seg{\bar\rr}{\bar\lll}})=\sum_{0\leq \ms_1\leq \rr_1-\lll_1}\cdots
\sum_{0\leq \ms_{N-1}\leq \rr_{N-1}-\lll_{N-1}}
\ \  \prod_{1\leq a\leq N-1}  \frac{1}{(\ms_a)!(\rr_a-\lll_a-\ms_a)!}\times \\ [10mm]
\ds \times\   \tSym_{\ \bar t_{\seg{\bar\rr}{\bar\lll}}}
\left(Z_{\bar\ms}({\bar t}_{\seg{\bar\rr}{\bar\lll}})
\ds\  \Pfm\sk{\F(\bar t_{\meg{\bar\lll}{\bar\lll+\bar\ms}})}\cdot
\Pfp\sk{\F(\bar t_{\meg{\bar\lll+\bar\ms}{\bar\rr}})}\right)
\end{array}
\end{equation}
where
\begin{equation}\label{Zser}
Z_{\bar\ms}({\bar t}_{\seg{\bar\rr}{\bar\lll}})=\prod_{a=1}^{N-2}\ \
\prod_{\substack{\lll_a+\ms_a<\ell\leq \rr_a \\
\lll_{a+1}<\ell'\leq \lll_{a+1}+\ms_{a+1}}} \frac{q^{-1}-q
t^{a+1}_{\ell'}/t^{a}_{\ell}}{1-t^{a+1}_{\ell'}/t^{a}_{\ell}}\,.
\end{equation}
\end{proposition}
To prove it one has to use the coproduct \r{gln-copr} and the property of the projection
\r{pr-prop}.

For a collections of variables $\bar t_{\meg{\bar\lll}{\bar\rr}}$ and the product
of the currents
$\tiF_i(t)$  \r{FFFl} the ordering is different.  Define a series
\begin{equation}\label{Zserm}
\hZ_{\bar\ss}({\bar t}_{\seg{\bar\rr}{\bar\lll}})=\prod_{a=1}^{N-2}\ \
\prod_{\substack{ \rr_a-\ss_a < \ell\leq \rr_a \\
  \lll_{a+1}< \ell' \leq \rr_{a+1}-\ss_{a+1}}} \frac{q-q^{-1}\
t^{a}_{\ell}\,/\,t^{a+1}_{\ell'}}{1-t^{a}_{\ell}\,/\,t^{a+1}_{\ell'}}\,.
\end{equation}
Then the ordered product of currents $\hF(\bar t_{\seg{\bar\rr}{\bar\lll}})$
 can be presented in the ordered
form (cf. \cite{KP-GLN})
\begin{equation}\label{dec-ff1m}
\begin{array}{c}
\ds \hF(\bar t_{\seg{\bar\rr}{\bar\lll}})=
\sum_{0\leq \ss_{N-1}\leq \rr_{N-1}-\lll_{N-1}}\cdots \sum_{0\leq \ss_1\leq \rr_1-\lll_1}
\ \  \prod_{1\leq a\leq N-1}  \frac{1}{(\ss_a)!(\rr_a-\lll_a-\ss_a)!}\times \\ [10mm]
\ds \times\   \tSym_{\ \bar t_{\seg{\bar\rr}{\bar\lll}}}
\left(\hZ_{\bar\ss}({\bar t}_{\seg{\bar\rr}{\bar\lll}})
\ds\  \hPfm\sk{\tiF(\bar t_{\seg{\bar\rr}{\bar\rr-\ss}})}\cdot
\hPfp\sk{\tiF(\bar t_{\seg{\bar\rr-\bar\ss}{\bar\lll}})}\right).
\end{array}
\end{equation}
Note that the rational series \r{Zser} and \r{Zserm}  when $\bar\ms+\bar\ss=\bar\rr-\bar\lll$
 are expansions of the same rational function
in the different zones. The set of the formal variables $\bar t_{\bar n}$
corresponds to the multisets \r{sord} and \r{hsord} in \r{dec-ff1} and \r{dec-ff1m}
respectively.

\subsection{Composed currents and strings}

Following the general strategy \cite{KP-GLN} and according to \cite{DKh} we  introduce the
{\it composed currents\/} associated to the currents
$\ff_{j,i}(t)$ for $i<j$.
The currents $\ff_{i+1,i}(t)$, $i=1,\ldots N-1$  coincide with
$\ff_i(t)$, cf.~\r{DF-iso}.
The coefficients of the series $\ff_{j,i}(t)$
belong to the completion $\overline U_F$ of the algebra $U_F$.
The composed currents $\ff_{j,i}(t)$ can be defined inductively using the formulas
\begin{equation}\label{rec-f1}
\ff_{j,i}(t)\, =\,
-\,\mathop{\rm res}\limits_{w=t}\ff_{a,i}(t)\ff_{j,a}(w)\,\frac{dw}{w}\, =\,
\mathop{\rm res}\limits_{w=t}\ff_{a,i}(w)\ff_{j,a}(t)\,\frac{dw}{w}
\end{equation}
for any $a=i+1,\ldots, j-1$. This relations are equivalent to
\begin{equation}\label{rec-f111}
\begin{split}
\ff_{j,i}(t)&=
\oint \ff_{a,i}(t) \ff_{j,a}(w)\ \frac{dw}{w}-
\oint \frac{q-q^{-1}t/w}{1-t/w}
\;\ff_{j,a}(w) \ff_{a,i}(t)\ \frac{dw}{w}\,,\\
\ff_{j,i}(t)&=
\oint \ff_{a,i}(w) \ff_{j,a}(t)\ \frac{dw}{w}-
\oint  \frac{q-q^{-1}w/t}{1-w/t}
\;\ff_{j,a}(t) \ff_{a,i}(w)\ \frac{dw}{w}\,.\\
\end{split}
\end{equation}
In \r{rec-f111}
$\oint \frac{dw}{w} g(w)=g_0$ for any formal series $g(w)=\sum_{n\in\ZZ}g_n z^{-n}$.
Using \r{gln-com} we can calculate
the residues in \r{rec-f1}:
\begin{equation}\label{res-in}
\ff_{j,i}(t)=(q-q^{-1})^{j-i-1}\ff_{j-1}(t) \ff_{j-2}(t)\cdots \ff_{i+1}(t)\ff_{i}(t)\,.
\end{equation}
Calculating formal integrals in \r{rec-f111} we obtain
expressions for the composed currents in the form
\begin{equation}\label{rec-f112}
\begin{split}
\ff_{j,i}(t)&=
\ff_{a,i}(t)\ff_{j,a}[0]-q\ff_{j,a}[0]\ff_{a,i}(t)-
(q-q^{-1})\sum_{k< 0}\ff_{j,a}[k]\,\ff_{a,i}(t)\,t^{-k}\,,\\
\ff_{j,i}(t)&=
\ff_{a,i}[0]\ff_{j,a}(t)-q^{-1}\ff_{j,a}(t)\ff_{a,i}[0]-
(q-q^{-1})\sum_{k\geq 0}\ff_{j,a}(t)\,\ff_{a,i}[k]\,t^{-k}\,,
\end{split}
\end{equation}
which are useful for the calculation of the projections of the composed currents.
By the similar procedure the composed currents $\tiF_{j,i}(t)$ associated
with the currents \r{DF2} were defined in \cite{KP-GLN}:
\begin{equation}\label{res-t}
\tiF_{j,i}(t)=(q-q^{-1})^{j-i-1}\tiF_{i}(t) \tiF_{i+1}(t)\cdots \tiF_{j-2}(t)\tiF_{j-1}(t)\,.
\end{equation}
The analytical properties of the products of the
composed currents, used in the paper, are presented in Appendix~\ref{anal-pr1}.
The analytical properties of the currents $\tiF_{j,i}(t)$ are given in the
Appendix A of the paper \cite{KP-GLN}.

Calculating the universal weight functions \r{uwf} and \r{uwf1} the special
products of the composed currents appear, which were called in \cite{KP}  as {\it strings}.
Firstly, we define the strings made from the currents $\ff_{j,i}(t)$.
Consider again a collection of segments $\seg{\bar\rr}{\bar\lll}$ and  associated
 set of variables $\bar t_{\seg{\bar\rr}{\bar\lll}}$.
Let $j=\textrm{min}(a)$ such that $\rr_b=\lll_b$ for $b=1,\ldots,a-1$.
Let $\bar\ms$ be a set of nonnegative integers, which satisfy the inequalities
\begin{equation}\label{restm}
0\leq\ms_a\leq \rr_a-\lll_a\,,\quad a=j+1,\ldots,N-1\,,
\end{equation}
and the non-increasing admissibility conditions
\begin{equation}\label{non-incm}
\rr_{j}-\lll_j=\ms_{j}\geq \ms_{j+1}\geq \ldots \geq \ms_{N-2}\geq \ms_{N-1}\geq \ms_N=0\,.
\end{equation}
Define an ordered product of the composed currents:
\begin{equation}\label{stringm}
\SF^j_{\bar\ms}({\bar t}^j_{\seg{\rr_j}{\lll_j}})=
\prod^{\longrightarrow}_{j+1\leq a\leq N}\sk{
\prod^{\longrightarrow}_{\rr_j-\ms_{a-1}< \ell\leq \rr_j-\ms_{a}}
\ff_{a,j}(t^j_\ell)}\,.
\end{equation}
The string \r{stringm} depends
 only on the variables of the type $j:\{ t_{\lll_j+1}^j,\ldots,t^j_{\rr_j}\}$,
 corresponding to the segment $\meg{\lll_j}{\rr_j}$.
The set of nonnegative integers $\bar\ms$ satisfying the admissibility condition
\r{non-incm} divides the segment $\meg{\lll_j}{\rr_j}$ into $N - j$ subsegments
$\meg{\rr_j-\ms_{a-1}}{\rr_j-\ms_a}$ for $a=j+1,\ldots,N$.
This division defines the product of the composed currents $\F^j_{\bar\ms}({\bar
t}^j_{\seg{\rr_j}{\lll_j}})$
in the string \r{stringm}.

Recall the construction of the strings from the composed currents $\tiF_{j,i}(t)$.
Consider again a collection of segments $\seg{\bar\rr}{\bar\lll}$ and  associated
 set of variables $\bar t_{\seg{\bar\rr}{\bar\lll}}$.
Let $j=\textrm{max}(a)$ such that $\rr_b=\lll_b$ for $b=a+1,\ldots,N-1$.
Let $\bar\ss$ be a set of nonnegative integers, which satisfy the inequalities
\begin{equation}\label{rest}
\rr_a-\lll_a\geq \ss_a\geq 0\,,\quad a=1,\ldots,j\,,
\end{equation}
and the non-increasing admissibility conditions
\begin{equation}\label{non-inc}
\rr_{j}-\lll_j=\ss_{j}\geq \ss_{j-1}\geq \ldots \geq \ss_{2}\geq \ss_1\geq \ss_0=0\,.
\end{equation}
Define an ordered product of the composed currents
\begin{equation}\label{string}
\ticSF^j_{\bar\ss}({\bar t}^j_{\seg{\rr_j}{\lll_j}})=
\prod^{\longleftarrow}_{j\geq a\geq 1}\sk{
\prod^{\longleftarrow}_{\lll_j+\ss_{a}\geq \ell> \lll_j+\ss_{a-1}}
\tiF_{j+1,a}(t^j_\ell)}
\end{equation}
for the string associated with the currents $\tiF_i(t)$.
The string \r{string} also depends
 only on the variables of the type $j:\{t^j_{\lll_j+1},\ldots, t_{\rr_j}^j\}$,
 corresponding to the segment $\seg{\lll_j}{\rr_j}$.
The set of nonnegative integers $\bar\ss$ satisfying the admissibility condition
\r{non-inc} divides the segment $\seg{\lll_j}{\rr_j}$ into $j$ subsegments
$\seg{\lll_j+\ss_a}{\lll_j+\ss_{a-1}}$ for $a=1,\ldots,j$.
This division defines the product of the composed currents in the string \r{string}.
  Projection of the string \r{stringm}
will be presented below. Projection of the string \r{string} was calculated in \cite{KP-GLN}.

For two sets of variable $\bar u$ and $\bar v$ we define the series
\begin{equation}\label{rat-Y}
\begin{split}
\ds U(u_1,\ldots,u_k;v_1,\ldots,v_k)&=\ds \prod_{1 \le i \le
k}\frac{v_i/u_i}{1-v_i/u_i} \;\;
\prod_{1 \le m < n \le k}\frac{q^{-1}-qv_m/u_{n}}{1-v_m/u_{n}}\,,\\
V(u_1,\ldots,u_k;v_1,\ldots,v_k)&=\ds \prod_{1 \le i \le
k}\frac{v_i/u_i}{1-v_i/u_i}\;\;
\prod_{1 \le n < m \le k}\frac{q-q^{-1}v_{n}/u_m}{1-v_{n}/u_m}\,.
\end{split}
\end{equation}
For the sets $\bar\ms$ and $\bar\ss$ which satisfy restrictions
 \r{restm}, \r{non-incm} and  \r{rest}, \r{non-inc}
  we define a series depending on the set of the
 variables  $\bar t_{\meg{\bar\lll}{\bar\lll+\bar\ms}}$ and  $\bar t_{\seg{\bar\rr}{\bar\rr-\bar\ss}}$:
\begin{equation}\label{rat-Xm}
Y(\bar t_{\seg{\bar\lll+\bar\ms}{\bar\lll}})=
\prod_{a=j+1}^{N-1}
U(t^{a-1}_{\lll_{a-1}+\ms_{a-1}-\ms_{a}+1},\ldots,t^{a-1}_{\lll_{a-1}+\ms_{a-1}};
t^{a}_{\lll_{a}+1},\ldots,t^{a}_{\lll_a+\ms_{a}} )
\end{equation}
and
\begin{equation}\label{rat-X}
X(\bar t_{\seg{\bar\rr}{\bar\rr-\bar\ss}})=
\prod_{a=1}^{j-1}
V(t^{a+1}_{\rr_{a+1}-s_{a+1}+s_{a}},\ldots,t^{a+1}_{\rr_{a+1}-s_{a+1}+1};
t^{a}_{\rr_{a}},\ldots,t^{a}_{\rr_{a}-s_{a}+1} )\,.
\end{equation}
When $j=N-1$ in \r{rat-Xm} we set $Y(\cdot)=1$ and if $j=1$ in \r{rat-X} we set $X(\cdot)=1$.

\subsection{Recurrence relation}\label{recur}
Let $\bar n$  be the set of nonnegative integers
$\bar n=\{n_{1},\ldots, n_{N-1}\}$.
We claim that the universal weight functions
$\calW^{N-1}(\bar{t}_{\segg{\bar{n}}})$ and $\ticalW^{N-1}(\bar{t}_{\segg{\bar{n}}})$
satisfies the recurrence relations
\begin{proposition}\label{prop45}${}$

\noindent
(i) The universal weight function {\rm \r{uwf1}} satisfies the recurrence relation
\begin{equation}\label{main-recm}
\begin{split}
\calW^{N}(\bar{t}_{\megg{\bar{n}}})&= \sum_{\bar\ms'}
\prod_{a=1}^{N-1}\frac{1}{(\ms_a-\ms_{a+1})!(n_a-\ms_{a})!}\times\\[5mm]
&\times\ \overline{\rm Sym}_{\ \bar t_{\megg{\bar n}}}
\sk{Z_{\bar\ms'}(\bar t_{\megg{\bar n'}})\cdot Y(\bar
t_{\megg{\bar\ms}})\cdot {\Pfp}\sk{\SF^{1}_{\bar\ms}(\bar
t^{1}_{\megg{n_{1}}})} \cdot
\calW^{N-1}(\bar{t}_{\meg{\bar\ms'}{\bar n'}})},
\end{split}
\end{equation}
where the sum is taken over all collections
$\bar\ms'=\{\ms_{2},\ms_{3},\ldots,\ms_{N-1}\}$ such
that $n_{1}=\ms_{1}\geq \ms_{2}\geq \cdots\geq \ms_{N-1}\geq \ms_N=0$
for $0\leq \ms_a\leq n_a$, $a=2,\ldots,N-1$ and
$\bar\ms=\{\ms_{1},\ms_{2},\ldots,\ms_{N-1}\}$, $\bar n'=\{n_2,\ldots,n_{N-1}\}$.

\noindent
(ii) \cite{KP-GLN}  The universal weight function {\rm\r{uwf}} satisfies the recurrence relation
\begin{equation}
\begin{split}
&\ticalW^{N}(\bar{t}_{\segg{\bar{n}}})= \sum_{\bar\ss'}
\prod_{a=1}^{N-1}\frac{1}{(\ss_a-\ss_{a-1})!(n_a-\ss_{a})!}\times\label{main-rec}\\[5mm]
&\times\ \overline{\rm Sym}_{\ \bar t_{\segg{\bar n}}}
\sk{\hZ_{\bar\ss'}(\bar t_{\segg{\bar n'}})\cdot X(\bar t_{\seg{\bar
n}{\bar n-\bar \ss}})\cdot {\hPfp}\sk{\ticSF^{N-1}_{\bar\ss}(\bar
t^{N-1}_{\segg{n_{N-1}}})} \cdot
\ticalW^{N-1}(\bar{t}_{\segg{\bar n'-\bar \ss'}})},
\end{split}
\end{equation}
where the sum is taken over all collections
$\bar \ss'=\{\ss_1,\ss_2,\ldots,\ss_{N-2}\}$, such
that $n_{N-1}=\ss_{N-1}\geq s_{N-2}\geq \cdots\geq s_1\geq s_0=0$
for $0\leq s_a\leq n_a$, $a=1,\ldots,N-2$ and
$\bar \ss=\{\ss_1,\ldots,\ss_{N-2},\ss_{N-1}\}$, $\bar n'=\{n_1,\ldots,n_{N-2}\}$.

\end{proposition}
{\it  Proof}\ \  of \r{main-rec} is
given in \cite{KP-GLN}.

We will sketch the proof of the recurrence
 relation \r{main-recm}.
We have a decomposition
$\F(\bar t_{\segg{\bar n}})=
\F(\bar t^{1}_{\segg{ n_{1}}})\F(\bar t_{\segg{\bar n'}})$,
where the first factor depends  only on the variables of the type $1$
and the second factor does not depend on the variables of this type:
\begin{equation*}
\begin{split}
\F(\bar t^{1}_{\segg{ n_{1}}})&=
\ff_{1}(t^{1}_{1})\cdots \ff_{1}(t^{1}_{n_{1}}),
\\
\F(\bar t_{\segg{\bar n'}})&=
\ff_{2}(t^{2}_{1})\cdots \ff_{2}(t^{2}_{n_{2}})\cdots
\ff_{N-1}(t^{N-1}_{1})\cdots
\ff_{N-1}(t^{N-1}_{n_{N-1}})\,.
\end{split}
\end{equation*}
We apply  to the product $\F(\bar t_{\segg{\bar n'}})$  the  ordering procedure of
Proposition~\ref{decomp} and substitute the result into \r{uwf}:
\begin{equation}\label{rr1}
\begin{split}
\mathcal{W}^{N}(\bar t_{\segg{\bar n}})&=
\sum_{\bar m'}\prod_{a=2}^{N-1}\frac{1}{\ms_a!(n_a-\ms_a)!}\times
\\
&\tSym_{\ \bar t_{\segg{\bar n'}}}
\left(Z_{\bar\ms'}(\bar t_{\segg{\bar n'}})
\ \Pfp \Big( \F(\bar t^{1}_{\segg{m_{1}}})
 \Pfm \sk{\F(\bar t_{\segg{\bar\ms'}})}\Big)
\cdot \mathcal{W}^{N-1}(\bar t_{\meg{\bar\ms'}{\bar n'}})\right).
\end{split}
\end{equation}
The sum is taken over all possible collections of the nonnegative integers
$\bar m' = \{\ms_{2},\ldots,$ $\ms_{N-1}\}$
 such that
 $\ms_a\leq n_a $, $a=2,\ldots,N-1$  and
$q$-symmetrization is performed over the set of the variables
$\bar t_{\segg{\bar n'}}$.

Now the proof of the Proposition~\ref{prop45} is based on the
following
\begin{lemma}\label{lemma43a}
For any collection of positive integers $\bar n=\{n_{1},\ldots,n_{N-1}\}$
and collections of the nonnegative integers ${\bar\ms}=
\{\ms_{1},\ms_{2},\ldots, \ms_{N-1}\}$,
 ${\bar\ms'}=
\{\ms_{2},\ldots,\ms_{N-1}\}$ such that all $\ms_a\leq n_a$ for
$a=2,\ldots,N-1$ and $\ms_{1}=n_{1}$ we have
\begin{equation}\label{rr2a}
\begin{split}
\Pfp&\sk{ \F(\bar t^{1}_{\segg{m_{1}}})
 \Pfm \sk{\F(\bar t_{\segg{\bar\ms'}})}}=
 \prod_{a=1}^{N-1}\frac{1}{(\ms_a-\ms_{a+1})!}\
\tSym_{\ \bar t_{\segg{\bar\ms}}}
\sk{Y(\bar t_{\segg{\bar\ms}})\cdot
\Pfp\sk{\SF^{1}_{\bar\ms}(\bar t^{1}_{\segg{m_{1}}})}}
\end{split}
 \end{equation}
 if $\ms_{1}\geq \ms_{2}\geq\cdots\geq \ms_{N-1}\ge m_N = 0$. Otherwise the right hand side of
 {\rm \rf{rr2a}} is equal to zero.
\end{lemma}
 The series $Y(\bar t_{\segg{\bar\ms}})$ is
defined by \rf{rat-Xm} as a particular case with $j=1$ and $\bar\lll=\bar0$.
\smallskip

{}From the definition \r{Zser} the series
$Z_{\bar\ms'}(\bar t_{\segg{\bar n'}})$ is symmetric with respect to permutations of the variables
of the same type from the set $t_{\segg{\bar\ms}}$ and the universal weight function
$\calW^{N-1}(\bar t_{\meg{\bar\ms'}{\bar n'}})$ does not depend on
the variables $t_{\segg{\bar\ms}}$. After substitution \r{rr2a} into \r{rr1}
 we can include the  series $Z_{\bar\ms'}(\bar t_{\segg{\bar n'}})$
and $\calW^{N-1}(\bar t_{\meg{\bar\ms'}{\bar n'}})$ in the
$q$-symmetrization $\overline{\rm Sym}_{\ \bar t_{\segg{\bar\ms}}}(\cdot)$
and using \r{exa3} replace the double $q$-symmetrization by a
single one:
\begin{equation*}
\tSym_{\ \bar t_{\segg{\bar n'}}}\ \tSym_{\
\bar t_{\segg{\bar\ms}}}(\ \cdot\ )=
\prod_{a=2}^{N-1}\ms_a!\ \tSym_{\ \bar t_{\segg{\bar n}}}(\ \cdot\ )\,.
\end{equation*}
The Proposition~\ref{prop45} is proved.
\hfill$\square$
\medskip

We shift the proof of the Lemma~\ref{lemma43a} to the Appendix~\ref{techB}.

\subsection{Iteration of the recurrence relation}
Let $\admis{\ms}=\{\ms_j^i\}$ and
$\admis{\ss}=\{\ss_i^j\}$ for $1\leq i\leq j\leq N-1$
be two collections of the nonnegative integers.
We say that collections  $\admis{\ms}$ and $\admis{\ss}$ are $\bar{n}$ admissible,
if they satisfy conditions of admissibility \r{admis-h} and \r{admis-m} respectively.
We also follow the convention $\ms^j_N=\ss_0^j=0$ for $j=1,...,N-1$.
Collections of integers $\admis{\ms}$ and
$\admis{\ss}$ can be visualized as triangular matrices
\begin{equation}\label{mad-mat}
\admis{\ms}=
\left(
\begin{array}{ccccc}
\ms_{1}^{1}& \ms_{2}^{1}&\ldots&\ms_{N-2}^{1}&\ms_{N-1}^{1}\\[2mm]
0&\ms_{2}^{2}&\ldots&\ms_{N-2}^{2}&\ms_{N-1}^{2}\\[2mm]
&&\ddots&\vdots&\vdots\\[2mm]
&0&&\ms^{N-2}_{N-2}&\ms^{N-2}_{N-1}\\[2mm]
&&&&\ms_{N-1}^{N-1}
\end{array}
\right)
\
\begin{array}{c}
0=\ms_N^{1}\!\!\leq \ms^{1}_{N-1}\leq \ldots\,\leq \ms_{1}^{1}\,,  \\[2mm]
0=\ms_N^{2}\!\!\leq \ms_{N-1}^{2}\leq \ldots \!\leq \ms_{2}^{2}\,, \\[2mm]
\vdots\\[2mm]
0=\ms^{N-2}_N\leq \ms^{N-2}_{N-1}\leq \ms^{N-2}_{N-2}\,,\\[2mm]
0=\ms^{N-1}_N\leq \ms_{N-1}^{N-1}\,,
\end{array}
\end{equation}
and
\begin{equation}\label{ad-mat}
\admis{s}=
\left(
\begin{array}{ccccc}
s_{1}^{1}&&&&\\[2mm]
s^{2}_1&s^{2}_{2}&&0&\\[2mm]
\vdots&\vdots&\ddots&&\\[2mm]
 s_{1}^{N-2}&s_{2}^{N-2}&\ldots&s_{N-2}^{N-2}&\\[2mm]
s_{1}^{N-1}& s_{2}^{N-1}&\ldots&s_{N-2}^{N-1}&s_{N-1}^{N-1}
\end{array}
\right)\
\begin{array}{c}
0=s^1_0\leq s_{1}^{1}\,, \\[2mm]
0=s^2_0\leq s^{2}_{1}\leq s^{2}_2\,,\\[2mm]
\vdots\\[2mm]
0=s_0^{N-2}\!\!\leq s_{1}^{N-2}\leq \ldots \!\leq s_{N-2}^{N-2}\,,   \\[2mm]
0=s_0^{N-1}\!\!\leq s_{1}^{N-1}\leq \ldots\,\leq s_{N-1}^{N-1}\,.
\end{array}
\end{equation}

Let $\bar\ms^j$ and $\bar\ss^j$, $j=1,\ldots,N-1$ be the $j$-th lines of the admissible matrices
  $\admis{\ms}$ and $\admis{\ss}$.
Define a collection of vectors
\begin{equation}\label{mi-def}
\bm^j=\bar\ms^{1}+\bar\ms^{2}+\cdots+\bar\ms^{j-1}+\bar\ms^{j}\,,\quad
 j=1,\ldots,N-1\,,
\end{equation}
\begin{equation}\label{si-def}
\bs^j=\bar\ss^{j}+\bar\ss^{j+1}+\cdots+\bar\ss^{N-2}+\bar\ss^{N-1}\,,\quad
 j=1,\ldots,N-1
\end{equation}
with non-negative integer components.
Set $\bm^0=\bar 0$ and $\bs^N=\bar 0$.
Denote by $\bfm^j_a$ and $\bfs^j_a$
 components of the vectors
$\bm^j$ and $\bs^j$:
\begin{equation*}
\bm^j=\{n_1,n_2,\ldots,n_j,\ms^{1}_{j+1}+\cdots+\ms^{j}_{j+1},\ldots,
\ms^{1}_{N-1}+\cdots+\ms^{j}_{N-1}\}\,,
\end{equation*}
\begin{equation*}
\bs^j=\{\ss^{j}_{1}+\cdots+\ss^{N-1}_{1},\ldots,\ss^{j}_{j-1}+\cdots+\ss^{N-1}_{j-1},
n_{j},\ldots,n_{N-2},n_{N-1}\}\,.
\end{equation*}
 Note that $\bfm^j_j-\bfm^{j-1}_{j}=\ms^{j}_{j}$ and
according to admissibility
condition \r{admis-h} and \r{admis-m} $\bm^{N-1}=\bs^1=\bar n$.

The iteration of the recurrence relations \r{main-recm} and \r{main-rec}
gives the following
\begin{theorem}\label{main-th}${}$

\noindent
(i)
The weight function {\rm \r{uwf1}} can
be presented as a total $q$-sym\-metri\-za\-tion  of the sum over all
 $\bar{n}$ admissible matrices
$\admis{\ms}$ of the ordered products of the projections of strings {\rm \r{stringm}}
with rational coefficients:
\begin{equation}\label{mWt1}
\begin{split}
&\ds \calW^{N}(\bar{t}_{\segg{\bar{n}}]})=  \tSym
_{\ \bar t_{\segg{\bar n}}}
\sum_{\admis{\ms}}\left(\prod_{b=1}^{N-1}\prod_{a=b}^{N-1}
\frac{1}{(\ms^b_a-\ms^b_{a+1})!}
\right.\times\\ &\quad \times\ds \left.
\prod_{j=1}^{N-3}
Z_{\bar\ms^j}(\bar t_{\meg{\bm^{j-1}}{\bar n}})
\prod_{j=1}^{N-2}
Y(\bar t_{\meg{\bm^{j-1}}{\bm^j}})
\prod_{1\leq j\leq N-1}^{\longrightarrow}
\Pfp\sk{\SF_{\bar\ms^j}^j (\bar t^j_{\meg{\bfm^{j-1}_{j}}{\bfm^j_j}})}
\right).
\end{split}
\end{equation}
(ii) \cite{KP-GLN}
The weight function {\rm \r{uwf}} can
be presented as a total $q$-sym\-metri\-za\-tion  of the sum over all
 $\bar{n}$ admissible matrices
$\admis{s}$ of the ordered products of the projections of strings {\rm \r{string}}
with rational coefficients:
\begin{equation}\label{Wt1}
\begin{split}
&\ds \ticalW^{N}(\bar{t}_{\segg{\bar{n}}]})=  \tSym
_{\ \bar t_{\segg{\bar n}}}
\sum_{\admis{\ss}}\left(\prod_{b=1}^{N-1}\prod_{a=1}^b
\frac{1}{(\ss^b_a-\ss^b_{a-1})!}
\right.\times\\ &\quad \times\ds \left.
\prod_{j=3}^{N-1}
\hZ_{\bar\ss^j}(\bar t_{\segg{\bar n-\bs^{j+1}}})
\prod_{j=2}^{N-1}
X(\bar t_{\seg{\bar n-\bs^{j+1}}{\bar n-\bs^{j}}})
\prod_{N-1\geq j\geq 1}^{\longleftarrow}
\hPfp\sk{\ticSF_{\bar\ss^j}^j (\bar t^j_{\segg{\ss^j_j}})}
\right).
\end{split}
\end{equation}
\end{theorem}

The rational series  in  \r{mWt1} and \r{Wt1}
may be gathered into a single multi-variable series. Define
\begin{equation}\label{gen-serm}
\begin{split}
\mathcal{Y}_{\admis{\ms}}&(\bar t_{\segg{\bar n}})=
\prod_{j=1}^{N-3}
Z_{\bar\ms^j}(\bar t_{\meg{\bm^{j-1}}{\bar n}})
\prod_{j=1}^{N-2}
Y(\bar t_{\meg{\bm^{j-1}}{\bm^j}})=\\
&=\prod_{a=2}^{N-1}\prod_{b=1}^{a-1}\prod_{\ell=0}^{\ms^b_a-1}
\frac{t^{a}_{\bfm^{b}_{a}-\ell}/t^{a-1}_{\bfm^{b}_{a-1}-\ell}}
{1-t^{a}_{\bfm^{b}_{a}-\ell}/t^{a-1}_{\bfm^{b}_{a-1}-\ell}}
\prod_{\ell'=\bfm^b_{a-1}-\ell+1}^{n_{a-1}}
\frac{q^{-1}-qt^{a}_{\bfm^{b}_{a}-\ell}/t^{a-1}_{\ell'}}
{1-t^{a}_{\bfm^{b}_{a}-\ell}/t^{a-1}_{\ell'}}
\end{split}
\end{equation}
and
\begin{equation}\label{gen-ser}
\begin{split}
\mathcal{X}_{\admis{\ss}}&(\bar t_{\segg{\bar n}})=
\prod_{j=3}^{N-1}
\hZ_{\bar\ss^j}(\bar t_{\segg{\bar n-\bs^{j+1}}})
\prod_{j=2}^{N-1}
X(\bar t_{\seg{\bar n-\bs^{j+1}}{\bar n-\bs^{j}}})=\\
&=\prod_{b=2}^{N-1}\prod_{a=1}^{b-1}\prod_{\ell=1}^{s^b_a}
\frac{t^{a}_{\ell+n_a-\bfs^{b}_{a}}/t^{a+1}_{\ell+n_{a+1}-\bfs^{b}_{a+1}}}
{1-t^{a}_{\ell+n_a-\bfs^{b}_{a}}/t^{a+1}_{\ell+n_{a+1}-\bfs^{b}_{a+1}}}
\prod_{\ell'=1}^{\ell+n_{a+1}-\bfs^b_{a+1}-1}
\frac{q-q^{-1}t^{a}_{\ell+n_a-\bfs^{b}_{a}}/t^{a+1}_{\ell'}}
{1-t^{a}_{\ell+n_a-\bfs^{b}_{a}}/t^{a+1}_{\ell'}}\,.
\end{split}
\end{equation}

We can formulate the following corollary of the Theorem~\ref{main-th}.
\begin{corollary}\label{corfromth1}
The weight functions {\rm \r{uwf1}}  and {\rm \r{uwf}} can be presented in the following
compact form:
\begin{equation}\label{mWt-fin}
\ds \calW^{N}(\bar{t}_{\segg{\bar{n}}})=  \tSym
_{\ \bar t_{\segg{\bar n}}}
\sum_{\admis{\ms}}\left(
\frac{\mathcal{Y}_{\admis{\ms}}(\bar t_{\segg{\bar n}})
}{\prod\limits_{a\geq b}(\ms^b_a-\ms^b_{a+1})!}
\prod^{\longrightarrow}_{1\leq j\leq N-1}
\Pfp\sk{\SF_{\bar\ms^j}^j (\bar t^j_{\meg{\bfm^{j-1}_{j}}{\bfm^j_j}})}
\right),
\end{equation}
\begin{equation}\label{Wt-fin}
\ds \ticalW^{N}(\bar{t}_{\segg{\bar{n}}})=  \tSym
_{\ \bar t_{\segg{\bar n}}}
\sum_{\admis{\ss}}\left(
\frac{\mathcal{X}_{\admis{s}}(\bar t_{\segg{\bar n}})
}{\prod\limits_{a\leq b}(s^b_a-s^b_{a-1})!}
\prod^{\longleftarrow}_{N-1\geq j\geq 1}
\hPfp\sk{\ticSF_{\bar s^j}^j (\bar t^j_{\segg{s^j_j}})    }
\right).
\end{equation}
\end{corollary}

Theorem~\ref{main-th} reduces the calculation of the universal weight function
to the calculation of the projections of the strings.

\subsection{Projection of composed currents and strings}

In this subsection we will formulate several statements about projections of the
strings. Their proofs can be found in papers \cite{KP,KP-GLN}.
First, define two types of the screening operators
\begin{equation}\label{scren}
S_A\ B=B\cdot A- q A\cdot B\,,\quad \hS_A\ B=B\cdot A- q^{-1} A\cdot B\,.
\end{equation}
We use these operators when $A$ are zero modes and $B$ are total currents
$\ff_i(t)$ or $\tiF_i(t)$. In this case these operators can be related
to the standard coproduct \r{coprL} via associated adjoint action (cf. details in \cite{KP-GLN}).

First relation in \r{rec-f112} for the currents $\ff_i(t)$ and analogous relation for the
currents $\tiF_i(t)$ yield
\begin{equation}\label{pro-in}
\begin{split}
\Pfp\sk{\ff_{j,i}(t)}&=
S_{\ff_{j,j-1}[0]}\bigl(\Pfp(\ff_{j-1,i}(t))\bigr)\,,\\
\hPfp\sk{\tiF_{j,i}(t)}&=
\hS_{\tiF_{i+1,i}[0]}\bigl(\hPfp(\tiF_{j,i+1}(t))\bigr)\,,\quad i<j-1\,.
\end{split}
\end{equation}
Iterating these relations we may write the projections for the composed currents
in terms of the successive application of the screening operators
\begin{equation*}
\begin{split}
\Pfp&\sk{\ff_{j+1,i}(u)}=S_{j}
                       S_{{j-1}}
                       \cdots
                       S_{{i+1}} \sk{ \Pfp( {\ff}_{i}(u)) }\\
                       &=
                       S_{j}
                       S_{{j-1}}
                       \cdots
                       S_{{i+1}}\left(   u\sk{u^{-1}{\ff}_{i}(u)} ^{(+)}\right)=
                       u\left(u^{-1}S_{j}
                       S_{{j-1}}
                       \cdots
                       S_{{i+1}} \sk{  {\ff}_{i}(u) }\right)^{(+)}\,,
                       \end{split}
\end{equation*}
\begin{equation*}
\begin{split}
\hPfp&\sk{\tiF_{j+1,i}(u)}=\hS_{i}
                       \hS_{{i+1}}
                       \cdots
                       \hS_{{j-1}} \sk{ \hPfp( {\tiF}_{j}(u)) }\\
                       &=
                       \hS_{i}
                       \hS_{{i+1}}
                       \cdots
                       \hS_{{j-1}}\left(   u\sk{u^{-1}{\tiF}_{j}(u)} ^{(+)}\right)=
                       u\left(u^{-1}\hS_{i}
                       \hS_{{i+1}}
                       \cdots
                       \hS_{{j-1}} \sk{  {\tiF}_{j}(u) }\right)^{(+)}\,,
                       \end{split}
\end{equation*}
where we introduced shorthand notations $S_i\equiv S_{F_i[0]}$ and $\hS_i\equiv \hS_{\tiF_i[0]}$.

\medskip

For a set $\{u_1,\ldots,u_n\}$ of  formal variables we introduce
two sets of the rational functions
\begin{equation}\label{rat-f}
\begin{split}
\varphi_{u_m}(u;u_1,\ldots,u_n)&=\prod_{k=1,\ k\neq m}^{n}\frac{u}{u_m}\
\frac{u-u_k}{u_m-u_k}\prod_{k=1}^{n}\frac{q^{-1}u_m-qu_k}{q^{-1}u-qu_k}\,,\\
\hat\varphi_{u_m}(u;u_1,\ldots,u_n)&=\prod_{k=1,\ k\neq m}^{n}\frac{u}{u_m}\
\frac{u-u_k}{u_m-u_k}\prod_{k=1}^{n}\frac{qu_m-q^{-1}u_k}{qu-q^{-1}u_k}\,,
\end{split}
\end{equation}
which satisfy the normalization conditions
$\varphi_{u_m}(u_s;u_1,\ldots,u_n)=\hat\varphi_{u_m}(u_s;u_1,\ldots,u_n)=\delta_{ms}$. We set
\begin{equation}\label{long-cur}
\begin{split}
\ff_{j+1,i}(u;u_{1},\ldots,u_n)&=\ff_{j+1,i}(u)-
\sum_{m=1}^n \varphi_{u_m}(u;u_{1},\ldots,u_n)\ff_{j+1,i}(u_m)\,,\\
\tiF_{j+1,i}(u;u_{1},\ldots,u_n)&=\tiF_{j+1,i}(u)-
\sum_{m=1}^n \hat\varphi_{u_m}(u;u_{1},\ldots,u_n)\tiF_{j+1,i}(u_m)\,,
\end{split}
\end{equation}
for $1\leq i\leq j<N$ and analogously
\begin{equation}\label{cc-pr1}
\begin{split}
\Pfp\sk{\ff_{j+1,i}(u;u_{1},\ldots,u_n)}&=S_{j}
                       S_{j-1}
                       \cdots
                       S_{i+1}\sk{F_{i}(u)^{(+)}}-\\
-&\sum_{m=1}^n \varphi_{u_m}(u;u_{1},\ldots,u_n)S_j
S_{j-1}\cdots                       S_{i+1}\sk{F_{i}(u_m)^{(+)}}\ ,\\
\hPfp\sk{\tiF_{j+1,i}(u;u_{1},\ldots,u_n)}&=\hS_{i}
                       \hS_{i+1}
                       \cdots
                       \hS_{j-1}\sk{\tiF_{j}(u)^{(+)}}-\\
-&\sum_{m=1}^n \hat\varphi_{u_m}(u;u_{1},\ldots,u_n)\hS_i
\hS_{i+1}\cdots                       \hS_{j-1}\sk{\tiF_{j}(u_m)^{(+)}}\ .
\end{split}
\end{equation}

In the same way, as it is described in the Appendix~B of the paper~\cite{KP-GLN} one
can write down following formulas for the projections of the strings

\begin{equation}\label{pr-st1m}
\begin{split}
&\ds \Pfp\sk{\SF^j_{\bar\ms}(\bar t^j_{\seg{\rr_j}{\lll_j}})}=
\prod^{\longleftarrow}_{N\geq a\geq j+1}\sk{
\prod^{\longleftarrow}_{\rr_j-\ms_{a}\geq \ell> \rr_j-\ms_{a-1}}
\Pfp\sk{F_{a,j}(t^j_\ell;t^j_{\ell+1},\ldots,t^j_{\rr_j})}}\times\\
\\ &\ds\qquad\qquad
\prod_{\lll_j<\ell<\ell'\leq \rr_j}
\frac{q^{-1}-qt^j_{\ell}/t^j_{\ell'}}{1-t^j_{\ell}/t^j_{\ell'}}
\prod_{j+1\leq a\leq N}\sk{
\prod_{\rr_j-\ms_{a-1}<\ell<\ell'\leq \rr_j-\ms_a}
\frac{1-t^j_{\ell}/t^j_{\ell'}}{q-q^{-1}t^j_{\ell}/t^j_{\ell'}} }\,,
\end{split}
\end{equation}
\begin{equation}\label{pr-st1}
\begin{split}
&\ds \hPfp\sk{\ticSF^j_{\bar\ss}(\bar t^j_{\seg{\rr_j}{\lll_j}})}=
\prod^{\longrightarrow}_{1\leq a\leq j}\sk{
\prod^{\longrightarrow}_{\lll_j+\ss_{a-1}< \ell\leq \lll_j+\ss_{a}}
\hPfp\sk{\tiF_{j+1,a}(t^j_\ell;t^j_{\lll_j+1},\ldots,t^j_{\ell-1})}}\times
\\ &\ds\qquad\qquad
\prod_{\lll_j<\ell<\ell'\leq \rr_j}
\frac{q-q^{-1}t^j_{\ell'}/t^j_{\ell}}{1-t^j_{\ell'}/t^j_{\ell}}
\prod_{1\leq a\leq j}\sk{
\prod_{\lll_j+s_{a-1}<\ell<\ell'\leq \lll_j+s_a}
\frac{1-t^j_{\ell'}/t^j_{\ell}}{q^{-1}-qt^j_{\ell'}/t^j_{\ell}} }\,.
\end{split}
\end{equation}

Commutation relations between total and Cartan currents together with the
formulas \r{pr-st1m}, \r{pr-st1} imply the relations
which are necessary for the next section:
\begin{equation}\label{gc8m}
\begin{split}
&\ds \Pfp\sk{\SF^j_{\bar s}(\bar t^j_{\seg{\rr_j}{\lll_j}})}
 \prod_{\ell=\lll_j+1}^{\rr_j} k^+_{j}(t^j_\ell)
=  \prod_{j+1\leq a\leq N}\sk{
\prod_{\rr_j-\ms_{a-1}<\ell<\ell'\leq \rr_j-\ms_a}
\frac{1-t^j_{\ell}/t^j_{\ell'}}{q-q^{-1}t^j_{\ell}/t^j_{\ell'}}
}\times\\
&\ds\quad\qquad\times \prod^{\longleftarrow}_{N\geq a\geq j+1}\sk{
\prod^{\longleftarrow}_{\rr_j-\ms_{a}\geq \ell> \rr_j-\ms_{a-1}}
\Pfp\sk{F_{a,j}(t^j_\ell)}k^+_{j}(t^j_\ell)}\,,
\end{split}
\end{equation}
\begin{equation}\label{gc8}
\begin{split}
&\ds \hPfp\sk{\ticSF^j_{\bar\ss}(\bar t^j_{\seg{\rr_j}{\lll_j}})}
 \prod_{\ell=\lll_j+1}^{\rr_j} \tk^+_{j+1}(t^j_\ell)
=  \prod_{1\leq a\leq j}\sk{
\prod_{\lll_j+\ss_{a-1}<\ell<\ell'\leq \lll_j+\ss_a}
\frac{1-t^j_{\ell'}/t^j_{\ell}}{q^{-1}-qt^j_{\ell'}/t^j_{\ell}}
}\times\\
&\ds\quad\qquad\times \prod^{\longrightarrow}_{1\leq a\leq j}\sk{
\prod^{\longrightarrow}_{\lll_j+\ss_{a-1}< \ell\leq \lll_j+\ss_{a}}
\hPfp\sk{\tiF_{j+1,a}(t^j_\ell)}\tk^+_{j+1}(t^j_\ell)}\,.
\end{split}
\end{equation}

\section{Universal weight functions and $\LL$-operators}
\label{sec5}

In this section we use factorization formulas \r{gc8m} and \r{gc8} to present
off-shell Bethe vectors  \r{b-v-sym} and \r{b-v-sym1} in terms of $\LL$-operator's
entries. First we relate projections of the currents to the Gauss coordinates of $\LL$-operators.
It can be proved that for $i<j$
\begin{equation}\label{ide-pr}
\begin{split}
\Pfp\sk{F_{j,i}(t)}\,&=\,(q^{-1}-q)^{j-i-1}\ \FF^{+}_{j,i}(t)\,,\\
\hPfp\sk{\tiF_{j,i}(t)}\,&=\,(q-q^{-1})^{j-i-1}\ \tFF^{+}_{j,i}(t)\,.
\end{split}
\end{equation}
For $j=i+1$ this is Ding-Frenkel \cite{DF} relations. For $i<j-1$, formulas
\r{ide-pr} follow from the  relations
between $\LL$-operators
\begin{equation*}
\begin{split}
(q^{-1}-q)\ \FF^+_{j,i}(t)\,&=\, S_{j-1} \sk{\FF^{+}_{j-1,i}(t)}\,,\\
(q-q^{-1})\ \tFF^+_{j,i}(t)\,&=\, \hS_{i} \sk{\tFF^{+}_{j,i+1}(t)}
\end{split}
\end{equation*}
written in terms of Gauss coordinates (see details in \cite{KPT}).

\begin{lemma}\label{lemmaideal}
{}For any $c=1,\ldots,N$ denote by $I_c$ and $\hI_c$ the left ideals of
$U_q(\widehat{\mathfrak{b}}^+)$,
generated by the modes of $\EE^+_{i,j}(u)$ with $ i<j\leq c$
and by the modes of $\tEE^+_{i,j}(u)$ with $ c \leq i < j$. We have inclusions $0=I_1\subset
 I_{2}\subset\cdots\subset I_N$ and $0=\hI_N\subset
 \hI_{N-1}\subset\cdots\subset \hI_1$.
\begin{itemize}
\item[{\it (i)}] Consider $a$ and $b$ such that $a<b$, then we have
equalities:
\begin{equation*}
\begin{split}
\LL^+_{a,b}(t)\,&\equiv\, \FF^+_{b,a}(t)k^+_a(t)\, \mod \ I_c\,,\
\qquad\LL^+_{a,a}(t)\,\equiv\, k^+_a(t)\, \mod \
I_c\,,\quad a\leq c\,,\\
\LL^+_{a,b}(t)\,&\equiv\, \tFF^+_{b,a}(t)\tk^+_b(t)\, \mod \ \hI_c\,,\
\qquad \LL^+_{b,b}(t)\,\equiv\, \tk^+_b(t)\, \mod \
\hI_c\,,\quad b\geq c\,.
\end{split}
\end{equation*}
\item[{\it (ii)}] The left ideal $I_c$ is generated by  modes of
$\ \LL^+_{j,i}(u)$ with $i < j \leq c$;
the left ideal $\hI_c$ is generated by  modes of
$\ \LL^+_{j,i}(u)$ with $c \leq i < j$, and $I_N = \hI_1$.
\item[{\it (iii)}] For any  $a\geq c$ the modes of
$\LL^+_{c,a}(t)$ normalize the ideal $I_c$:
\begin{equation*}
I_c\cdot  \LL^+_{c,a}(t) \subset I_c\,.
\end{equation*}
For any  $a\leq c$ the modes of
$\LL^+_{a,c}(t)$ normalize the ideal $\hI_c$:
\begin{equation*}
\hI_c\cdot  \LL^+_{a,c}(t) \subset \hI_c\,.
\end{equation*}
\end{itemize}
\end{lemma}

\noindent
{\it Proof.}\ \ The items {\it (i)}\/ and {\it (ii)}\/ are proved
using the Gauss decompositions
\r{GF2}--\r{GE2} and \r{GF1}--\r{GE1}. Item  {\it (iii)}\/ follows
from {\it (ii)}\/ and $RLL$-relations \r{L-op-com}.
\hfill$\square$
\medskip
\begin{theorem}\label{th-Wt3} For any $\Uqgln$ module $V$ with a weight singular
vector $v$ we have the following formulas for modified weight functions
\begin{equation}\label{mWt3}
\begin{split}
\ds &{\bf w}_V^{N}(\bar{t}_{\segg{\bar{n}}})=\ds
\qsym(\bar t_{\segg{\bar n}})\
\tSym
_{\ \bar t_{\segg{\bar n}}}
\sum_{\admis{\ms}}\Bigg(
\frac{(q^{-1}-q)^{\sum\limits_{b=1}^{N-1}(n_b-\ms_b^b)}}
{ \prod\limits_{a\geq b}^{N-1}(\ms^b_a-\ms^b_{a+1})!}
\ \mathcal{Y}_{\admis{\ms}}(\bar t_{\segg{\bar n}})\times
  \\
 &\prod_{1\leq b\leq N-1}^{\longrightarrow}\Bigg(
\prod_{N-1\geq a\geq b}^{\longleftarrow}
\prod_{\ell= n_b-\ms^b_a+1}^{n_b-\ms^b_{a+1}}
\bigg( \LL^{+}_{b,a+1}(t^b_\ell)
\prod_{\ell'=n_b-\ms^b_a+1}^{\ell-1}
\frac{t^b_{\ell'}-t^b_{\ell}}{q^{-1}t^b_{\ell'}-qt^b_{\ell}} \bigg)\Bigg)
 \prod_{\ell=1}^{n_b-\ms^b_b}\LL^+_{b,b}(t^b_\ell)\Bigg)v
\end{split}
\end{equation}
and
\begin{equation}\label{Wt3}
\begin{split}
&\ds \hat{\bf w}_V^{N}(\bar{t}_{\segg{\bar{n}}})=\ds
\qsym(\bar t_{\segg{\bar n}})\
 \tSym_{\ \bar t_{\segg{\bar n}}}
\sum_{\admis{\ss}}\left(
\frac{(q-q^{-1})^{\sum\limits_{b=1}^{N-1}(n_b-\ss_b^b)}}{\prod\limits_{a\leq b}^{N-1}
(\ss^b_a-\ss^b_{a-1})!}
\ \mathcal{X}_{\admis{\ss}}(\bar t_{\segg{\bar n}})\times
\right.\\
 &\ds \prod\limits_{N-1\geq b\geq 1}^{\longleftarrow}\left.\!\!\Bigg(
\prod_{1\leq a\leq b}^{\longrightarrow}
\prod_{\ell=s^b_{a-1}+1}^{s^b_a}\bigg(\LL^{+}_{a,b+1}(t^b_\ell)
\prod_{\ell'=\ell+1}^{s^b_a}
\frac{t^b_\ell-t^b_{\ell'}}{q^{-1}t^b_\ell-qt^b_{\ell'}} \bigg)\Bigg)
\prod_{\ell=s^b_b+1}^{n_b}\LL^{+}_{b+1,b+1}(t^b_\ell)
\right)v
\end{split}
\end{equation}
(the ordering in the products over $\ell$ is not important, because of commutativity
of the entries of $\LL$-operators with equal matrix indices).
\end{theorem}
\smallskip
{\it Proof.}\ Let $V$ be $\Uqgln$-module with a weight singular vector $v$. Using
the relations \r{gc8m}, \r{gc8}, \r{ide-pr} and
the corollary~\ref{corfromth1}
we can write modified weight functions ${\bf
w}_V^{N}(\bar{t}_{\segg{\bar{n}}})$ and $\hat{\bf w}_V^{N}(\bar{t}_{\segg{\bar{n}}})$ in the following form:
\begin{equation}\label{mWt2}
\begin{split}
\ds &{\bf w}_V^{N}(\bar{t}_{\segg{\bar{n}}})=\ds
\qsym(\bar t_{\segg{\bar n}})\
\tSym
_{\ \bar t_{\segg{\bar n}}}
\sum_{\admis{\ms}}\left(
\frac{(q^{-1}-q)^{\sum_{b=1}^{N-1}(n_b-\ms_b^b)}}
{ \prod\limits_{a\geq b}^{N-1}(\ms^b_a-\ms^b_{a+1})!}
\ \mathcal{Y}_{\admis{\ms}}(\bar t_{\segg{\bar n}})\!\!\!\!
\prod_{1\leq b\leq N-1}^{\longrightarrow}\left(
\prod_{N-1\geq a\geq b}^{\longleftarrow}
  \right.\right.\\
 &\left.\left.
\times \prod\limits_{n_b-\ms^b_{a+1}\geq\ell> n_b-\ms^b_a}^{\longleftarrow}
\left( \FF^{+}_{a+1,b}(t^b_\ell)k^+_{b}(t^b_\ell)
\prod_{\ell'=n_b-\ms^b_a+1}^{\ell-1}
\frac{t^b_{\ell'}-t^b_{\ell}}{q^{-1}t^b_{\ell'}-qt^b_{\ell}} \right)\right)
 \prod_{\ell=1}^{n_b-\ms^b_b}k^+_{b}(t^b_\ell)\right)v
\end{split}
\end{equation}
and
\begin{equation}\label{Wt2}
\begin{split}
\ds &\hat{\bf w}_V^{N}(\bar{t}_{\segg{\bar{n}}})=\ds
\qsym(\bar t_{\segg{\bar n}})\
\tSym
_{\ \bar t_{\segg{\bar n}}}
\sum_{\admis{\ss}}\left(
\frac{(q-q^{-1})^{\sum_{b=1}^{N-1}(n_b-\ss_b^b)}}
{ \prod\limits_{a\leq b}^{N-1}(\ss^b_a-\ss^b_{a-1})!}
\ \mathcal{X}_{\admis{\ss}}(\bar t_{\segg{\bar n}})\!\!\!\!
\prod_{N-1\geq b\geq 1}^{\longleftarrow}\left(
\prod_{1\leq a\leq b}^{\longrightarrow}
  \right.\right.\\
 &\left.\left.
\times \prod_{\ss^b_{a-1}<\ell \leq \ss^b_a}^{\longrightarrow}
\left( \tFF^{+}_{b+1,a}(t^b_\ell)\tk^+_{b+1}(t^b_\ell)
\prod_{\ell'=\ell+1}^{\ss^b_a}
\frac{t^b_\ell-t^b_{\ell'}}{q^{-1}t^b_\ell-qt^b_{\ell'}} \right)\right)
 \prod_{\ell=\ss^b_b+1}^{n_b}\tk^+_{b+1}(t^b_\ell)\right)v\,,
\end{split}
\end{equation}
where series $\mathcal{Y}_{\admis{s}}(\bar t_{\segg{\bar n}})$ and
$\mathcal{X}_{\admis{s}}(\bar t_{\segg{\bar n}})$ are given by \r{gen-serm} and \r{gen-ser}
respectively. Then the statement of the Theorem follows from
\r{mWt2}, \r{Wt2} and Lemma~\ref{lemmaideal} in the same way as in the
paper \cite{KP-GLN}.\qed

As in \cite{KP-GLN}
we can simplify  formulas \r{mWt3} and \r{Wt3} by using
the  $q$-sym\-met\-ri\-za\-tion \r{qs1r} (see \cite{TV3}).
Denoting
$\ppp^b_a=n_a-\bfs^b_a=\ss^{a}_a+\cdots+\ss^{b-1}_a$ we formulate the
following corollary of the Theorem~\ref{th-Wt3}
\begin{corollary}\label{cor53}
The off-shell Bethe vectors for the quantum affine algebra $\Uqgln$ can be written as
\begin{equation}\label{mWt4}
\begin{split}
\ds {\bf w}_V^{N}(\bar{t}_{\segg{\bar{n}}})&=\ds  \qSym
_{\ \bar t_{\segg{\bar n}}}
\sum_{\admis{\ms}}\left((q^{-1}-q)^{\sum\limits_{a=1}^{N-1}(n_a-\ms_a^a)}
\prod_{a\leq b}^{N-1}
\frac{1}{[\ms^a_b-\ms^a_{b+1}]_q!}\  \ \times \right.\\
&\times\prod_{a=2}^{N-1}\prod_{b=1}^{a-1}\prod_{\ell=0}^{\ms^b_a-1}
\frac{t^{a}_{\bfm^{b}_{a}-\ell}/t^{a-1}_{\bfm^{b}_{a-1}-\ell}}
{1-t^{a}_{\bfm^{b}_{a}-\ell}/t^{a-1}_{\bfm^{b}_{a-1}-\ell}}
\ \ \prod_{\ell'=\bfm^b_{a-1}-\ell+1}^{n_{a-1}}
\frac{q^{-1}-qt^{a}_{\bfm^{b}_{a}-\ell}/t^{a-1}_{\ell'}}
{1-t^{a}_{\bfm^{b}_{a}-\ell}/t^{a-1}_{\ell'}}
\\
&\times\ds
\left.\prod_{1\leq a\leq N-1}^{\longrightarrow}\left(
\prod_{N-1\geq b\geq a}^{\longleftarrow}
\prod_{\ell= n_a-\ms^a_b+1}^{n_a-\ms^a_{b+1}}
\LL^{+}_{a,b+1}(t^a_\ell)\right)
 \prod_{\ell=1}^{n_a-\ms^a_a}\LL^+_{a,a}(t^a_\ell)\right)v
\end{split}
\end{equation}
and
\begin{equation}\label{Wt4}
\begin{split}
\ds \hat{\bf w}_V^{N}(\bar{t}_{\segg{\bar{n}}})&=\ds  \qSym
_{\ \bar t_{\segg{\bar n}}}
\sum_{\admis{\ss}}\left((q-q^{-1})^{\sum\limits_{b=1}^{N-1}(n_b-\ss_b^b)}
\prod_{a\leq b}^{N-1}
\frac{1}{[\ss^b_a-\ss^b_{a-1}]_q!}\  \ \times \right.\\
&\times\
\prod_{b=2}^{N-1}\prod_{a=1}^{b-1}\prod_{\ell=1}^{\ss^b_a}
\frac{1}{1-t^{a}_{\ell+\ppp^{b}_{a}}/t^{a+1}_{\ell+\ppp^{b}_{a+1}}}
\prod_{\ell'=1}^{\ell+\ppp^b_{a+1}-1}
\frac{q-q^{-1}t^{a}_{\ell+\ppp^{b}_{a}}/t^{a+1}_{\ell'}}
{1-t^{a}_{\ell+\ppp^{b}_{a}}/t^{a+1}_{\ell'}}
\\
&\times\ds \left.
\prod_{N-1\geq b\geq 1}^{\longleftarrow} \sk{\prod_{1\leq a\leq b}^{\longrightarrow}
\sk{\prod_{\ell=\ss^b_{a-1}+1}^{\ss^b_a} \LL^+_{a,b+1}(t^b_\ell) }
\prod_{\ell=\ss^b_b+1}^{n_b}\LL^+_{b+1,b+1}(t^b_\ell)}
\right)v\,.
\end{split}
\end{equation}
\end{corollary}
{\it Proof.}\ \
We prove this corollary only for \eqref{mWt4}, since \eqref{Wt4}
can be proved in the same way. Consider right hand side of formula \eqref{mWt3}.
Inside the total $q$-symmetrization
there is a symmetric series in the sets of variables
$\{t^{b}_\ell\}$ for $n_b-\ms^b_{a}+1\leq \ell\leq n_b-\ms^b_{a+1}$
for $b=1,\ldots,N-1$ and $a=b,\ldots,N-1$. It follows from the
commutativity of matrix elements of $\LL$-operators
$[\LL^+_{a,b}(t),\LL^+_{a,b}(t')]=0$ and the explicit form of the
series~\r{gen-serm}. Product $\prod\limits_{\ell'=n_b-\ms^b_a+1}^{\ell-1}
\frac{t^b_{\ell'}-t^b_{\ell}}{q^{-1}t^b_{\ell'}-qt^b_{\ell}}$ is
an inverse to the function $\beta(t^{b}_{[n_b - m^b_a, n_b - m^b_{a+1}]})$
defined in \eqref{betadef}. Consequently applying formulas \eqref{po-sim},
\eqref{relat} and using the explicit form of the  series \eqref{gen-serm}
 we obtain \eqref{mWt4}.\qed

\section{Connection between two weight functions}\label{sec6}

Let us calculate the image of the off-shell Bethe vector \r{mWt4} for the evaluation
homomorphism \r{eval2}. Let $M_\Lambda$ be a $\Uqqq$-module generated by a vector $v$,
satisfying the conditions $\ele_{a,a}\,v=q^{\weig_a}\,v$ and
$\ele_{a,b}\,v=0$ for $a<b$. Then $v$ is a singular weight vector of
the $\Uqgln$ evaluation module $M_{\lambda}(z)$. The action of the matrix elements of
$\LL$-operators in this module is given by the formulas
\begin{equation}\label{eval-act}
\begin{split}
\mathcal{E}v_z\sk{{\LL}_{a,b+1}^+(t)}&=(q-q^{-1})\ele_{b+1,a}\ele_{b+1,b+1}\equiv
(q-q^{-1})\check\ele_{b+1,a}\,,\\
\mathcal{E}v_z\sk{{\LL}_{a,a}^+(t)}v&=\sk{q^{\weig_a}-q^{-\weig_a}\frac{z}{t}}v=\lambda_a(t)v\,.
\end{split}
\end{equation}

\begin{proposition}\label{dontknowdontcare}
For any evaluation $\Uqgln$ module $M_{\lambda}(z)$ with singular weight vector $v$ we have
$$
\displaystyle {\bf w}_{M_{\lambda}(z)}(\bar{t}_{\segg{\bar{n}}}) =
{\mathbb{B}}_{M_{\lambda}(z)}(\bar t_{\segg{\bar n}})\,.
$$
\end{proposition}
{\it Proof.}\ Substituting \r{eval-act} into \r{mWt4} and using reordering of the factors
\begin{equation}\label{reorder}
\prod_{a=2}^{N-1}\prod_{\ell=1}^{n_a-\ms^a_a}\lambda_{a}(t^a_\ell)=
\prod_{a=2}^{N-1}\prod_{b=1}^{a-1}\prod_{\ell=0}^{\ms^b_a-1}\lambda_{a}(t^a_{\bfm^b_a-\ell})
\end{equation}
we obtain
\begin{equation}\label{mmWt4}
\begin{split}
\ds &{\bf w}_{M_{\lambda}(z)}^{N}(\bar{t}_{\segg{\bar{n}}})=\ds (q-q^{-1})^{\sum\limits_{a=1}^{N-1}n_a}
\sum_{\admis{\ms}}\left( \left( \prod_{1\leq b\leq N-1}^{\longrightarrow}
\prod_{N-1\geq a\geq b}^{\longleftarrow}
\frac{\check\ele_{a+1,b}^{\ms^b_a-\ms^b_{a+1}}}{[\ms^b_a-\ms^b_{a+1}]_q!} \right)v \right.\\
&\left.\times\ \qSym _{\ \bar t_{\segg{\bar n}}}
\prod_{a=2}^{N-1}\prod_{b=1}^{a-1}\prod_{\ell=0}^{\ms^b_a-1}
\frac{q^{\weig_a}t^{a}_{\bfm^{b}_{a}-\ell} - q^{-\weig_a} z  }
{t^{a-1}_{\bfm^{b}_{a-1}-\ell}- t^{a}_{\bfm^{b}_{a}-\ell}}
\ \ \prod_{\ell'=\bfm^b_{a-1}-\ell+1}^{n_{a-1}}
\frac{qt^{a}_{\bfm^{b}_{a}-\ell} -q^{-1}t^{a-1}_{\ell'}}
{t^{a}_{\bfm^{b}_{a}-\ell}-t^{a-1}_{\ell'}}
\right)\,.
\end{split}
\end{equation}

Reordering generators $\check\ele_{a+1,b}$ using Serre relations \r{ser-fin} we
literally obtain equation~\eqref{Wt555y}. \hfill $\Box$

\medskip

Denote by $J$ the left ideal of $\Uqgln$ defined in the
item {\it (ii)}\/ of the Lemma~{\ref{lemmaideal}}:
$J = I_N = \hat I_1$. This ideal is generated by the modes of the entries of
the $\LL$-operator $\LL^+_{j,i}(t)$ with $1 \le i < j \le N$.  We use the statement proved
in the Theorem 3 of \cite{KP-GLN}, which we formulate as the following lemma.
\begin{lemma}\label{KPT3}
Consider $\Uqbpp{}$ as an algebra over $\CC[[(q - 1)]]$, and two elements $\mathcal{A}$,
 $\mathcal{B}\in\Uqbpp{}$. If for any singular weight vector $v$ we have that $\mathcal{A}\ v
 = \mathcal{B}\ v$,
 then
$$
\mathcal{A} \equiv \mathcal{B}\hbox{ {\rm mod} }J.
$$
\end{lemma}
\medskip
\begin{theorem}\label{mainmain} Universal weight functions are subject to the following relations:
\begin{itemize}
\item[{\it (i)}] For each irreducible finite-dimensional $\Uqgln$-module $V$ with a
singular vector $v$ two weight functions are equal:
$$
{\bf w}_V^{N}(\bar t_{\segg{\bar n}})={\mathbb{B}}_V(\bar t_{\segg{\bar n}}).
$$
\item[{\it (ii)}]
$${\bf w}_V^{N}(\bar t_{\segg{\bar n}}) = \hat{\bf w}_V^{N}(\bar
t_{\segg{\bar n}}).$$
\item[{\it (iii)}] Consider $\Uqbpp{}$ as an algebra over $\CC[[(q - 1)]]$, then
$$
\qsym(\bar t_{\segg{\bar n}})
\ \calW^{N}(\bar{t}_{\segg{\bar{n}}})
\prod_{a=1}^{N-1}\prod_{\ell=1}^{n_{a}}k^+_{a}(t^a_\ell) \equiv {\mathbb{B}}_V(\bar t_{\segg{\bar
n}})\quad \hbox{ {\rm mod} } J.
$$
\item[{\it (iv)}]
$$\calW^{N}(\bar{t}_{\segg{\bar{n}}})
\prod\limits_{a=1}^{N-1}\prod\limits_{\ell=1}^{n_{a}}k^+_{a}(t^a_\ell)
\equiv\ticalW^{N}(\bar{t}_{\segg{\bar{n}}})\prod_{a=1}^{N-1}\prod_{\ell=1}^{n_{a}}\tk^+_{a+1}(t^{a}_\ell)
\quad\hbox{ {\rm mod} } J.$$
\end{itemize}
\end{theorem}
{\it Proof.}\ \
($i$) Since both ${\bf w}_V^{N}(\bar t_{\segg{\bar n}})$ and
${\mathbb{B}}_V(\bar t_{\segg{\bar n}})$ have the same comultiplication
properties (cf. \cite{KPT}), it follows from the
Proposition~\ref{dontknowdontcare} that for any set of
evaluation $\Uqgln$ modules $M_{\lambda_i}(z_i)$
$$
{\bf w}_{1_{g(z_0)}\otimes M_{\lambda_1}(z_1)\otimes\dots \otimes M_{\lambda_{k}(z_k)}}^{N}(\bar
t_{\segg{\bar n}}) ={\mathbb{B}}_{1_{g(z_0)}\otimes M_{\lambda_1}(z_1)\otimes\dots \otimes
M_{\lambda_{k}(z_k)}}(\bar t_{\segg{\bar n}})
$$
for any tensor product of evaluation modules and a one-dimensional
module $1_{g(z_0)}$, in which every $L_{ii}(z)$ acts by
multiplication on $g(z_0)$ (both weight functions in consideration
are trivial for one-dimensional modules).
Then we can apply classical result of~\cite{CP}: every irreducible
finite-dimensional $\Uqgln$-module with a singular vector $v$ is isomorphic to a
subquotient of a tensor product of one-dimensional modules and of evaluation
modules, generated by the tensor product of their weight singular vectors. The
singular vector corresponds to the image tensor product of singular vectors
within this isomorphism.

{\it (ii)}\/ It follows from {\it (i)}\/ and the fact proved
in \cite{KP-GLN} that
$$
\hat{\bf w}_V^{N}(\bar
t_{\segg{\bar n}}) = {\mathbb{B}}_V(\bar t_{\segg{\bar n}}).
$$

{\it (iii)}\/ It follows from {\it (i)}\/ and Lemma~\ref{KPT3}.

{\it (iv)}\/ It follows from {\it (ii)}\/ and Lemma~\ref{KPT3}.
\qed

\section*{Acknowledgement}

The authors thank S.Khoroshkin for numerous useful discussions.
This work was partially done when the second
author (S.P.) visited Laboratoire d'An\-necy-Le-Vieux de Physique Th\'eorique in 2006 and 2007.
These visits were possible due to the financial support of
the CNRS-Russia exchange program on mathematical physics.
He thanks LAPTH for the hospitality and stimulating scientific atmosphere.
His work was supported in part by RFBR grant 05-01-01086 and grant for support
of scientific schools NSh-8065.2006.2.

This paper is a part of PhD thesis
of A.S. which he is prepared in co-direction of S.P. and V.R. in
the Bogoliubov  Laboratory of Theoretical Physics, JINR, Dubna and in LAREMA,
D\'epartement de Math\'ematics, Universit\'e d'Angers. He is
grateful to the CNRS-Russia exchange program on mathematical physics  and  personally to J.-M. Maillet
for financial and general support of this thesis project.

\setcounter{section}{0}
\setcounter{subsection}{0}

\renewcommand{\thesection}{\Alph{section}}


\section{Analytical properties of composed currents}\label{anal-pr1}

In this appendix we reformulate the Serre relations in terms of the composed currents.
We start from the following properties of the `regularized' products of usual total currents:

\begin{itemize}
 \item[(i)] $B_1(z,w)=(q^{-1}z-qw)F_i(z)F_i(w)$ vanish at $z=w$: $B_1(z,z)=0$;
 \item[(ii)] $B_2(z_1,z_2,z_3)=(z_1-z_2)(qz_2-q^{-1}z_3)(q^{-1}z_1-qz_3)F_i(z_1)F_{i+1}
 (z_2)F_i(z_3)$ vanish at $z_1=z_2=q^{-2}z_3$: $B_2(z,z,q^2z)=0$;
 \item[(iii)] $B_3(z_1,z_2,z_3)=(qz_1-q^{-1}z_2)(z_2-z_3)(q^{-1}z_1-qz_3)F_{i+1}(z_1)F_i
 (z_2)F_{i+1}(z_3)$ vanish at $z_1=z_2=q^{2}z_3$: $B_3(z,z,q^{-2}z)=0$.
\end{itemize}

The property~(i) follows from the commutation relations~(3.10) written as $$B_1(z,w)=-B_1(w,z).$$
The properties~(ii) and (iii) can be obtained from the Serre relations~(3.11) and
delta-function identities~\cite{E}. Note that the antisymmetry of $B_1(z,w)$
implies $B_2(z,w,u)=-B_2(u,w,z)$ and $B_3(z,w,u)=-B_3(u,w,z)$.

Let us demonstrate how one can derive the commutation relation
between $F_i(z)$ and the composed current $F_{i+2,i}(w)$
defined as $F_{i+2,i}(w)=u^{-1}(u-w)F_i(u)F_{i+1}(w)\Big|_{u=w}=
(q^{-1}-q)F_{i+1}(w)F_i(w)$. Using $B_2(z,w,u)=-B_2(u,w,z)$ we
have\footnote{Rational function $\frac1{1-x}$ should be always understood
 as formal series $\sum_{n\geq0}x^n$.}
\begin{equation*}
 F_i(z)F_{i+2,i}(w)=w^{-1}\frac{z^{-1}}{(1-w/z)}\
 \frac{qz^{-1}}{(1-q^2w/z)}\ B_2(w,w,z)\,.
\end{equation*}
On the other hand
\begin{equation*}
\frac{q^{-1}-qz/w}{1-z/w}\ F_{i+2,i}(w)F_i(z)=w^{-1}
\frac{q^{-1}w^{-1}}{(1-q^{-2}z/w)}\ \frac{w^{-1}}{(1-z/w)}\ B_2(w,w,z)\,.
\end{equation*}
Now, the property~(ii) means
 $$\frac{qz^{-1}}{1-q^2w/z}\ B_2(w,w,z)=-\frac{q^{-1}w^{-1}}{1-q^{-2}z/w}\ B_2(w,w,z)\,.$$
  The equality
 $B_2(z,z,z)=0$, which follows from~(i), means
 $$\frac{z^{-1}}{1-w/z}\ B_2(w,w,z)=-
 \frac{w^{-1}}{1-z/w}\ B_2(w,w,z)\,.$$
This proves the commutation relation:
\begin{align}
 F_{i+1,i}(z)F_{i+2,i}(w)=\frac{q^{-1}-qz/w}{1-z/w}\
 F_{i+2,i}(w)F_{i+1,i}(z)\,. \label{FFi2}
\end{align}
It follows from this commutation relation, the commutation relation between $F_{i+1,i}(z)$ and
$F_{i,i-1}(w)$ and the definition $F_{i+2,i-1}(w)=(q^{-1}-q)F_{i+2,i}(w)
F_{i,i-1}(w)$ that
\begin{align}
 F_{i+1,i}(z)F_{i+2,i-1}(w)=F_{i+2,i-1}(w)F_{i+1,i}(z)\,. \label{FFi3}
\end{align}
Generalizing  relations~\eqref{FFi2}, \eqref{FFi3} we obtain:
\begin{proposition}
 For any $i>j>k>l$ we have the following commutation relations
\begin{align}
 F_{jk}(z)F_{ik}(w)&=\frac{q^{-1}-qz/w}{1-z/w}\ F_{ik}(w)F_{jk}(z)\,, \label{FFi2g} \\
F_{ik}(z)F_{ij}(w)&=\frac{q^{-1}-qz/w}{1-z/w}\  F_{ij}(w)F_{ik}(z)\,, \label{FFi22g} \\
 F_{jk}(z)F_{il}(w)&=F_{il}(w)F_{jk}(z)\,, \label{FFi3g} \\
\frac{q^{-1}-qw/z}{1-w/z}\ F_{ij}(z)F_{ij}(w)&=\frac{q^{-1}-qz/w}{1-z/w}\ F_{ij}(w)F_{ij}(z)\,.
\label{FFijij}
\end{align}
\end{proposition}

\section{Proof of the Lemma~\ref{lemma43a}}
\label{techB}

Before proving this  lemma  we formulate several preliminary
propositions.

For any $j=1,...,N-1$
denote by $U_{j}$ the subalgebra of $\U_f$ formed by the modes of
$\ff_{j}(t),\ldots,\ff_{N-1}(t)$. Let $U_j^\varepsilon=U_j\cap
{\mathrm Ker}\, \varepsilon$ be the corresponding augmentation
ideal. Let $\bms{j}=\{\ms_{j},\ms_{j+1},\ldots,\ms_{N-1}\}$
 be a collection of non-negative  integers satisfying
admissibility conditions: $\ms_{j}\geq \ms_{j+1}\geq \cdots \geq \ms_{N-1}\geq \ms_N=0$.
 Set $\bms{j+1}=\{\ms_{j+1},\ldots,\ms_{N-1}\}$.

\begin{proposition}\label{main-fact}
For the product  $\F(\bar t^{j}_{\segg{\ms_{j}}})$
and  string $\SF^{j+1}_{\bms{j+1}}({\bar t}^{j+1}_{\segg{\ms_{j+1}}})$
we have a formal series equality
\begin{equation}\label{ca16}
\begin{split}
\F(\bar t^{j}_{\segg{\ms_{j}}}&) \cdot\Pfm\sk{
\SF^{j+1}_{\bms{j+1}}({\bar
t}^{j+1}_{\segg{\ms_{j+1}}})} =
\prod_{i=j}^{N-1}\frac{1}{(\ms_{i}-\ms_{i+1})!}\times
\\
&\tSym_{\ t^{j}_{\segg{\ms_{j}}}}\Big(
\SF^{j}_{\bms{j}}({\bar t}^{j}_{\segg{\ms_{j}}})\ \prod_{a=j+1}^{N-1}
\tSym_{\ \bar t^{j+1}_{\meg{\ms_{j+1}-\ms_{a}}{\ms_{j+1}-\ms_{a+1}}}}
\\
&\qquad\qquad U(t^{j}_{{m_j-m_{j+1}+1}},\ldots,
t^{j}_{{m_j}};t^{j+1}_{1},\ldots,t^{j+1}_{m_{j+1}})\Big)\quad {\rm mod}\
\Pfm\sk{U^\coun_{j+1}}\cdot U_{j}\ ,
\end{split}
\end{equation}
where a symbol $\prod\limits_{a=j+1}^{N-1}\tSym_{\ \bar
t^{j+1}_{\meg{\ms_{j+1}-\ms_{a}}{\ms_{j+1}-\ms_{a+1}}}}$ is a composition of the
corresponding q-sym\-met\-ri\-za\-tions: $\;\;\tSym_{\ \bar
t^{j+1}_{\meg{0}{\ms_{j+1}-\ms_{j+2}}}} \cdots\; \tSym_{\ \bar
t^{j+1}_{\meg{m_{j+1}-m_{N-1}}{\ms_{j+1}}}}$. If admissibility conditions
$\ms_{j}\geq \ms_{j+1}\geq \cdots \geq \ms_{N-1}\geq0$ are not satisfied then
the right hand side of {\rm \r{ca16}} is zero modulo $\Pfm\sk{U^\coun_{j+1}}\cdot U_{j}$.
\end{proposition}
{\it Proof.}\  The Proposition~\ref{main-fact} can be proved along the same lines
as it was done in the Appendix~C of the paper~\cite{KP-GLN}.
To perform this proof the reader should use formulas for the commutations of
the composed currents gathered in the Appendix~\ref{anal-pr1} of
this paper.\qed

\begin{lemma}\label{qsym-pr2}
Rational series\ \  $\tSym_{\ \bar t^{j+2}_{\segg{\ms_{j+2}}}}\cdots
\tSym_{\ \bar t^{N-2}_{\segg{\ms_{N-2}}}}
\tSym_{\ \bar t^{N-1}_{\segg{\ms_{N-1}}}} Y(\bar t_{\bms{j+1}})$ is symmetric in each
group of variables $\{t^{j+1}_{\ms_{j+1}-\ms_a+1},\ldots,t^{j+1}_{\ms_{j+1}-\ms_{a+1}}\}$ for
$a=j+1,\ldots,N-1$.
\end{lemma}
{\it Proof.}\  This Lemma can be proved by induction. The case of $j = N-3$
follows from the fact proved in \cite{KP} that the $q$-symmetrization of the formal series
\r{rat-Y}
$\tSym_{\ \bar v}\ U(\bar u;\bar v)$ is the symmetric series with respect to the set
of variables $\bar u$. For $j < N-3$ induction follows from the same fact,
property of $q$-symmetrization \eqref{sym_div} and
the formula
\begin{align}
U(u_1,\dots,u_k;v_1,\dots,v_k) &= \Sigma(u_{s+1},\dots,u_k;v_1,\dots,v_s)\times\nn\\
&U(u_1,\dots,u_s;v_1,\dots,v_s) U(u_{s+1},\dots,u_k;v_{s+1},\dots,v_k)\,,\notag
\end{align}
where $0 \le s \le k$ and $\Sigma(u_{s+1},\dots,u_k; v_1,\dots, v_s)$ is a symmetric
function with respect to each set of variables $\{u_{s+1},\dots,u_k\}$ and $\{v_1, \dots,
v_s\}$.\qed

\begin{proposition}\label{fact3}
There is a formal series equality
\begin{equation}\label{pr-ap1}
\begin{split}
&\F(\bar t^{j}_{\segg{\ms_{j}}}) \cdot
\tSym_{\ {\bar t}_{\segg{\bms{j+1}}}}
\sk{Y(\bar t_{\segg{\bms{j+1}}})\cdot\Pfm\sk{
\SF^{j+1}_{\bms{j+1}}({\bar t}^{j+1}_{\segg{\ms_{j+1}}})}}
=\\
&\qquad=\frac{1}{(\ms_{j}-\ms_{j+1})!}\
\tSym_{\ \bar t_{\segg{\bms{j}}}}
\sk{Y(\bar t_{\segg{\bms{j}}})\cdot
\SF^{j}_{\bms{j}}({\bar t}^{j}_{\segg{\ms_{j}}})}\quad
{\rm mod}\ \Pfm\sk{U^\coun_{j+1}}\cdot U_{j}
\end{split}
\end{equation}
is admissibility conditions
$\ms_{j}\geq \ms_{j+1}\geq \cdots \geq \ms_{N-1}\geq0$ are satisfied; otherwise
the right hand side of the equality {\rm \r{pr-ap1}} is zero modulo $\Pfm\sk{U^\coun_{j+1}}\cdot U_{j}$.
\end{proposition}
{\it Proof.}\  Lemma~\ref{qsym-pr2} implies the following relation
\begin{align}\label{f439}
&\tSym_{\ \bar t^{j+1}_{\segg{\ms_{j+1}}}}\bigg(
\tSym_{\ \bar t^{j+2}_{\segg{\ms_{j+2}}}}\cdots
\tSym_{\ \bar t^{N-1}_{\segg{\ms_{N-1}}}} Y(\bar t_{\bms{j+1}})\times\nn\\
&\quad\times
\prod_{a=j+1}^{N-1}
\tSym_{\ \bar t^{j+1}_{\meg{\ms_{j+1}-\ms_{a}}{\ms_{j+1}-\ms_{a+1}}}}
U\big(\bar t^{j}_{\seg{\ms_j-\ms_{j+1}}{\ms_j}}; \bar
t^{j+1}_{\segg{\ms_{j+1}}}\big)\bigg)=\\
&\qquad =\prod_{a=j+1}^{N-1}(\ms_a-\ms_{a+1})!\ \
\tSym_{\ \bar t_{\segg{\bms{j+1}}}} {Y(\bar t_{\bms{j}})}\,.\nn
\end{align}
The $q$-symmetrization $\tSym_{\ {\bar t}_{\segg{\bms{j+1}}}}$ in the left hand side of
\r{pr-ap1}
do not affect the formal series depending on the variables $\bar t^{j}_{\segg{\ms_{j}}}$.
 This means that the left hand side of \r{pr-ap1} can be written in the form
$\tSym_{\ {\bar t}_{\segg{\bms{j+1}}}}
\bigg({Y(\bar t_{\segg{\bms{j+1}}})\cdot\F(\bar t^{j}_{\segg{\ms_{j}}}) \cdot\Pfm\sk{
\SF^{j+1}_{\bms{j+1}}({\bar t}^{j+1}_{\segg{\ms_{j+1}}})}}\bigg)$. Then,
substituting~\eqref{ca16} into this expression, and using~\eqref{f439} we obtain
the right hand side of \r{pr-ap1}.\hfill$\Box$

\begin{proposition}\label{pr-rr6}
The projection $\Pfm \sk{\F(\bar t_{\segg{\bms{2}}})}$ can be presented as
\begin{equation}\label{rr6}
\prod\limits_{a=2}^{N-1}\frac{1}{(\ms_{a}-\ms_{a+1})!} \tSym_{\ \bar t_{\segg{\bms{2}}}}
\bigg(Y(\bar t_{\segg{\bms{2}}})
\
\ \Pfm\Big({\SF^{2}_{\bms{2}}(\bar t^2_{\segg{\ms_{2}}})  }\Big)\bigg)
  \mod\  \Pfm\sk{ U_{3}^\varepsilon}\cdot U_{2}
\end{equation}
 if admissibility conditions
$\ms_{2}\ge\ms_{3}\ge\cdots\ge\ms_{N-2}\ge\ms_{N-1} \ge\ms_N=0$
are satisfied and is zero modulo $\Pfm\sk{ U_{3}^\varepsilon}\cdot
U_{2}$ otherwise.
\end{proposition}
{\it Proof.}\  The statement of the Proposition~\ref{pr-rr6} is a particular
case of the following formal series equality
\begin{align}\label{rr7}
\Pfm \sk{\F(\bar t_{\segg{\bms{2}}})} &=
\prod\limits_{a=j+1}^{N-1}\frac{1}{(\ms_{a}-\ms_{a+1})!}\Pfm\bigg(\F(\bar
t^2_{\segg{\ms_{2}}})\cdots \F(\bar t^{j-1}_{\segg{\ms_{{j-1}}}})\times\\
&\F(\bar t^j_{\segg{\ms_{j}}})\
\tSym_{\ \bar t_{\segg{\bms{j+1}}}}
\Big(Y(\bar t_{\segg{\bms{j+1}}})
\ {\SF^{j+1}_{\bms{j+1}}(\bar t^{j+1}_{\segg{\ms_{j+1}}})  }\Big)\bigg)
  \mod\ \Pfm\sk{
U_{j+2}^\varepsilon}\cdot U_{2}\,.\nn
\end{align}
We prove \r{rr7} by induction. Indeed it is correct for $j = N-1$. Suppose
that it is valid for some $j \le N - 1$. Using the property~\eqref{pr-prop} and
taking into account~\eqref{pgln} one can replace the string
$\SF^{j+1}_{\bms{j+1}}(\bar t^{j+1}_{\segg{\ms_{j+1}}})$ in the right hand
side of~\eqref{rr7} by its negative projection  $\Pfm\big(\SF^{j+1}_{\bms{j+1}}(\bar
t^{j+1}_{\segg{\ms_{j+1}}})\big)$. From the Proposition~\ref{fact3} we
obtain the right hand side of~\eqref{rr7} for $j-1$, using also the facts that
$\Pfm\Big(\F(\bar
t^2_{\segg{\ms_{2}}})\cdots \F(\bar t^{j-1}_{\segg{\ms_{{j-1}}}})
\Pfm(U_{j+1}^\varepsilon)\cdot U_j\Big) \subset \Pfm(U_{j+1}^\varepsilon)\cdot U_2$ and
$\Pfm(U_{j+2}^\varepsilon)\cdot U_2 \subset \Pfm(U_{j+1}^\varepsilon)\cdot U_2$.
If the admissibility condition is not satisfied for some $j\geq2$, namely, $\ms_{j}<\ms_{j+1}$,
then according to the Proposition~\ref{fact3} the product
$\F(\bar t^j_{\segg{\ms_{j}}})\
\tSym_{\ \bar t_{\segg{\bms{j+1}}}}
\Big(Y(\bar t_{\segg{\bms{j+1}}})
\ {\SF^{j+1}_{\bms{j+1}}(\bar t^{j+1}_{\segg{\ms_{j+1}}})  }\Big)$
is zero modulo elements from $\Pfm(U_{j+1}^\varepsilon)\cdot U_j$.
Due to the commutativity of the modes of the currents $F_k[n]$, $k=2,\ldots,j-1$ with any element
from $U_{j+1}^\varepsilon$ we conclude that in this case the right hand side of \r{rr6} is zero
modulo $\Pfm(U_{j+1}^\varepsilon)\cdot U_2 \subset \Pfm(U_{3}^\varepsilon)\cdot U_2$.
\qed

\noindent{\it Proof of the Lemma~\ref{lemma43a}}.
For the calculation of projection $\Pfp \sk{
\F(\bar t^{1}_{\segg{n_{1}}})
 \Pfm \sk{\F(\bar t_{\segg{\bms{2}}})}}$ we substitute
 $\Pfm \sk{\F(\bar t_{\segg{\bms{2}}})}$
 using \r{rr6}. Since  any element of $U^\varepsilon_{3}$ commutes with
 the current  $\ff_{1}(t)$, elements of $\Pfm\sk{ U_{3}^\varepsilon}\cdot
U_{2}$ do not contribute: $\Pfp \Big(
\F(\bar t^{1}_{\segg{n_{1}}})
\cdot \Pfm\sk{ U_{3}^\varepsilon}\cdot
U_{2} \Big)=0$.
Then we have
\begin{equation*}
\prod_{a=2}^{N-1}\frac{1}{(\ms_{a}-\ms_{a+1})!}\
\Pfp \sk{
\F(\bar t^{1}_{\segg{n_{1}}})\
\tSym_{\ \bar t_{\segg{\bms{2}}}}
\sk{
Y(\bar t_{\segg{\bms{2}}})  \Pfm\!\sk{\SF^{2}_{\bms{2}}(\bar
t^2_{\segg{\ms_{2}}})  }}}\!,
\end{equation*}
where $\ms_N = 0$.
The latter expression
is  non-zero due to the Proposition \ref{fact3} only if $\ms_{2}\leq n_{1}$ and is equal to the
right hand side of \r{rr2a}. \qed

\end{document}